\documentclass[11pt,leqno,a4paper,english]{amsart}
\usepackage{amssymb,amsmath,amsfonts,amsthm,amscd}
\usepackage{leftidx}
\usepackage{tikz}
\usepackage[latin1]{inputenc}

\usepackage{mathrsfs}

\usepackage{hyperref}
\hypersetup{
    colorlinks=true,       
    linkcolor=blue,          
    citecolor=blue,        
    filecolor=magenta,      
    urlcolor=cyan           
}

\usepackage{subfigure,graphicx}
\usepackage{xcolor}
\usepackage{xspace}

\usepackage{enumerate}

\usepackage[all]{hypcap}

\let \Re \relax
\DeclareMathOperator{\Re}{Re}

\newcommand{\mb}[1]{\ensuremath{\mathbb{#1}}}

\newcommand{\R}{{\mb{R}}}
\newcommand{\C}{{\mb{C}}}

\newcommand{\eps}{\varepsilon}

\newcommand{\E}{\ensuremath{\mathcal E}}
\newcommand{\G}{\ensuremath{\mathcal G}}
\renewcommand{\H}{\ensuremath{\mathcal H}}

\renewcommand{\L}{\ensuremath{\mathcal L}}

\renewcommand{\d}{\ensuremath{\partial}}

\usepackage{subfigure,graphicx}
\usepackage{color}
\usepackage{xspace}

\newtheoremstyle{note}{} {}{\itshape}{-6pt}{\bf}{. --}{ }{}

\newtheorem{theorem}{Theorem}[section]
\newtheorem{proposition}[theorem]{Proposition}
\newtheorem{lemma}[theorem]{Lemma}
\newtheorem{corollary}[theorem]{Corollary}

\newtheorem*{theo*}{Theorem}

\newcounter{theorembiss}

\newtheorem{defi}[theorem]{Definition}
\newtheorem{rema}[theorem]{Remark}

\numberwithin{equation}{section}

\subjclass[2010]{35A18, 35Lxx, 35Q93, 93Dxx}

\title[Boundary Sidewise Observability of the Wave Equation] 
{Boundary Sidewise Observability of the Wave Equation} 

 \author{Belhassen Dehman}
 \address{Belhassen Dehman. D\'epartement de Math\'ematiques, Facult\'e
  des sciences de Tunis $\&$ Enit-Lamsin, Universit\'e de Tunis El Manar, 2092 El
  Manar, Tunisia. }
\email{belhassen.dehman@fst.utm.tn}

 \author{Enrique Zuazua}
 \address{Enrique Zuazua. [1] Chair for Dynamics, Control and Numerics - Alexander von Humboldt-Professorship, Department of Data Science, Friedrich-Alexander-Universit\"at Erlangen-N\"urnberg,
91058 Erlangen, Germany ,
\newline \indent\hskip 7pt
[2] Chair of Computational Mathematics, Fundaci\'{o}n Deusto,
48007 Bilbao, Basque Country, Spain,
\newline \indent\hskip 7pt
[3] Departamento de Matem\'{a}ticas,
Universidad Aut\'{o}noma de Madrid,
28049 Madrid, Spain.}
\email{enrique.zuazua@fau.de}
\date{\today}

\begin{document}
\begin{abstract} The wave equation on a bounded domain of $\R^{n}$ with non homogeneous boundary Dirichlet data or sources  supported on a subset of the boundary is considered. We analyze the problem of observing the source  out of boundary measurements done away from its support.

 We first show that observability inequalities may not hold unless an infinite number of derivatives are lost, due to the existence of solutions that are arbitrarily concentrated near the source. 
 
 We then establish observability inequalities in Sobolev norms, under a suitable microlocal geometric condition on the support of the source and the measurement set,   for sources  fulfilling  pseudo-differential conditions that exclude these concentration phenomena. 
 
 The proof relies on  microlocal arguments and is essentially based on the use of microlocal defect measures. \end{abstract} 

\maketitle
\tableofcontents
\section{Introduction}
\subsection{General setting}\label{setting}
Let $\Omega$ be a bounded open domain of $\R^{n}$ with boundary $\partial\Omega$ of class $\mathcal{C}^{\infty}$. We set 
$$ 
\L = \R\times \Omega \quad \text{and } \quad \d\L = \R\times \d\Omega .
$$ 
We also introduce  $A=(a _{ij}(x))$,  a $n \times n$ matrix  of  $\mathcal{C}^{\infty}$ coefficients, symmetric,  uniformly definite positive on a neighborhood of $\Omega$. 

Finally, we take $g  \in H^{1}(\d\L) $ and we assume  that  $g$ is compactly supported in time in the interval $(0,+\infty)$.

We consider then the following wave system 
\begin{equation}
\left\{ 
\begin{array}{c}
P_{A}u=\partial _{t}^{2}u-\sum_{i,j=1}^{n}\partial_{x_{j}}(a _{ij}(x)\partial_{x_{i}}u)=0\quad \text{in } \L 
\\ 
\\
u(t,.)=g(t,.)\quad \text{on } \d\L 
\\
\\
u(0,.)=\partial _{t}u(0,.)=0 \quad \text{in } \Omega.
\end{array}%
\right.  \label{waveequation}
\end{equation}
This system is well posed in the classical energy space
$
C^{0}(\R, H^{1}(\Omega)) \cap C^{1}(\R, L^{2}(\Omega))
$
equipped with the energy norm $\sup_{t\in \R}Eu(t)$, where $$Eu(t)= \Vert u(t,.)\Vert_{H^{1}(\Omega)}^{2} + \Vert \d_{t}u(t,.)\Vert_{L^{2}(\Omega)}^{2},$$and
$$
\Vert u(t,.)\Vert_{H^{1}(\Omega)}^{2} = \sum_{i,j=1}^n \int_\Omega a_{ij}(x) \partial_{x_{i}}u \partial_{x_{j}}u dx,
$$ 
see  \cite{La-Lio-Triggiani}. Actually, the solution $u$ vanishes  for  $t\leq 0$.

More precisely, the following energy estimate holds
\begin{equation}\label{energy}
\sup_{t\in \R}Eu(t) \leq C \vert\vert g\vert\vert_{H^{1}(\d\L )}^{2},
\end{equation}
together with the added hidden regularity property of the trace of the normal derivative
\begin{equation}\label{b-energy}
\Vert \d_{n}u_{\vert \partial \Omega}\Vert_{L^{2}((0,a)\times \partial \Omega )} \leq C_{a} \Vert g\Vert_{H^{1}(\d\L )},
\end{equation}
valid  for all $a >0$.
\begin{rema}
 The constant appearing  in estimate \eqref{energy} and \eqref{b-energy}  depend on the metric attached to $A=(a_{ij}(x))_{ij}$, on the geometry of the domain $\Omega$ and, for \eqref{b-energy}, also on  on the time-horizon $a>0$.
\end{rema}

\subsection{Geometry of the domain $\Omega$}
In this paper, we will deal with a particular class of domains $\Omega$. This fact is made precise  in the following condition.
\smallskip

\textbf{Assumption A1}

{\it We assume that there exists  a strictly concave (with respect to the metric attached to the matrix $A=(a_{ij}(x))_{ij}$) open non empty subset $O$ of the boundary $\d\Omega$, $\overline{O}\neq \d\Omega$.}
\smallskip

Geometrically, this guarantees that every geodesic of $\Omega$ that is tangent to $O$ at some point $m_{0}$, has an order of tangency equal to 1; locally near this point and except for $m_{0}$, this geodesic lives in $\Omega$.

\begin{figure}[h!]\label{geometryO}
 \begin{center}
 \subfigure
{\includegraphics[scale=0.9]{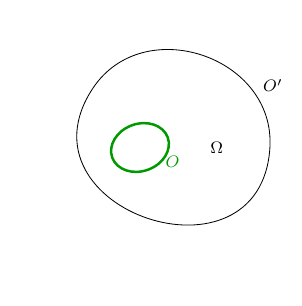}}
\quad\quad\quad\quad
\subfigure
{\includegraphics[scale=0.9]{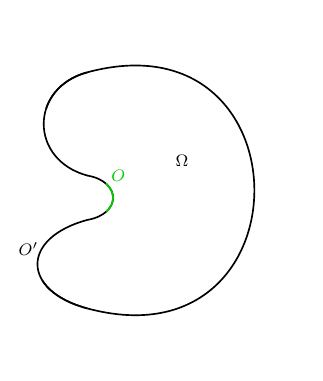}}
 \caption{Examples of strictly concave boundary subset $O$}
\end{center}
\end{figure}

For instance, if $A=Id$ , this simply says that there exists a neighborhood $V$ of $O$ in $\R^{n}$, such that  the set $V\setminus\Omega $ is strictly convex. See Fig.\ref{geometryO}.
\begin{rema}

\begin{enumerate}

\item Assumption A1, implicitly, substantially limits the class of domains $\Omega$ under consideration. For example, this condition excludes convex domains $\Omega$. Indeed, for subsets $O$ of the boundary of $\Omega$  to exist, so that they fulfil the assumption A1, the geometry of $\Omega$ needs to allow for some concavity zones of its boundary, as illustrated in Figure  \ref{geometryO}, and this excludes many domains $\Omega$.

\item In the literature, sets $O$  fulfilling assumption A1 are sometimes said to be diffractive with respect to the metric attached to $A=(a_{ij}(x))_{ij}$. 
\end{enumerate}
\end{rema}

\subsection{Motivation}
From now, we will work under assumption A1. Let then $O'$ be a non empty open subset of $\partial\Omega$ such that  $\overline{O}\cap \overline{O'} = \emptyset$. We set 
$$
\Gamma = \R\times O,  \quad\quad\Gamma' = \R\times O' , 
$$ 
and for $a >0$, 
$$
\L_{a}= (0,a)\times \Omega, \quad \Gamma_{a}= (0,a)\times O \quad \text{and} \quad \Gamma'_{a}= (0,a)\times O' .
$$
In addition, we assume throughout the whole paper that the boundary data $g$ is supported in $\overline{\Gamma}_{M}= [0,M]\times \overline{O}$ for some $M>0$.

The aim of this paper is to analyze whether it is possible to observe the boundary data or source $g$ in \eqref{waveequation} from measurements done on the normal derivative  $\d_{n}u_{\vert \Gamma'}$ on the subset $\Gamma'$ of the boundary.  In other words, we are seeking for an estimate of the type
\begin{equation}\label{obs}
\Vert g \Vert_{H^{1}(\Gamma_{M})} \leq C \Vert \d_n u_{\vert \d\Omega}\Vert_{L^{2}(\Gamma'_{a})} ,
\end{equation}
for some  $a \geq M$.
\begin{figure}[h!]
  \begin{center}
\includegraphics[scale=0.25]{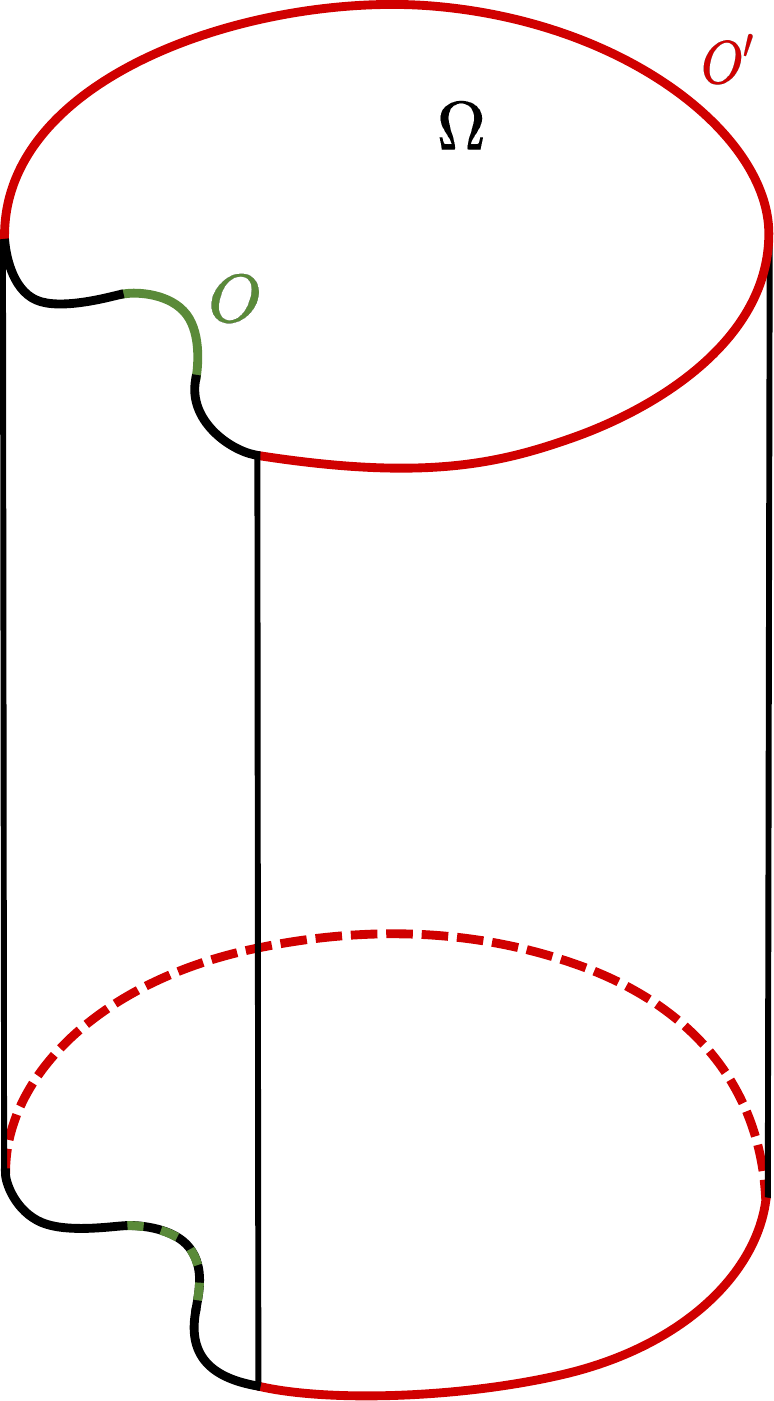}
   \caption{
   Cylindrical domain where waves evolve. In green the support of the source $g$ to be identified, and in  red the subset of the boundary where measurements are done.}
  \end{center}
\end{figure}

Estimate \eqref{obs} is the sidewise observability inequality object of analysis in this paper.

According to the Rellich inequality it is  well known that the right hand side term of \eqref{obs} is bounded above by 
$$
\Vert u \Vert^2_{a} =: \sup_{t\in [0,a]}Eu(t)=\sup_{t\in [0,M]}Eu(t)=\Vert u \Vert^2_{M}.
$$ 
More precisely, for every $a>0$, there exists $C_{a}>0$ such that  every solution $u$ of \eqref{waveequation} satisfies
\begin{equation}\label{continu}
\Vert \d_n u_{\vert \d\Omega}\Vert_{L^{2}(\Gamma'_{a})} \leq C_{a} \Vert u \Vert_{M}.
\end{equation}

Therefore, a necessary condition for an estimate of the form \eqref{obs} to hold is that the boundary data $g$ under consideration needs  to be observable out of the total interior energy $\Vert u \Vert_{M}$, namely, the existence of a constant $C>0$ such that  
\begin{equation}\label{obs-int}
\Vert g \Vert_{H^{1}(\Gamma_{M})} \leq C \Vert u\Vert_{M} .
\end{equation}

However, as we shall see, this inequality does not hold without additional structural conditions on the source term $g$ under consideration. Indeed, in Theorem \ref{Obs-Loss-1} and Corollary \ref{Obs-Loss-3}, we construct  sequences of invisible sources  $(g_{k})$ whose energy is essentially localized on the elliptic and/or glancing set of the boundary, such that
\begin{equation}\label{invisible}
\Vert g_{k} \Vert_{H^{1}(\Gamma_{M})}\rightarrow1, \quad g_{k} \rightharpoonup 0 \quad \text{in}  \,\,\, H^{1}, \quad \Vert u_{k} \Vert_{M} \longrightarrow 0,
\end{equation}
which, of course, are an impediment for \eqref{obs-int} to occur.

In fact, as we shall see, even the weaker version
\begin{equation}\label{obs-loss}
\Vert g \Vert_{H^{s}(\Gamma_{M})} \leq C \Vert \d_n u_{\vert \d\Omega}\Vert_{H^{1}(\Gamma'_{a})}
\end{equation}
may not for hold for any $s \leq 1$.

The lack of such sidewise observability inequalities is genuinely a multi-d phenomenon (see section \ref{sec 6}). 
By the contrary, as shown in  \cite{Sar-Zua} and \cite{Zuazua} by means of sidewise energy estimates, in 1-d , inequality \eqref{obs-int} holds for $BV$ coefficients and under natural conditions on  the length of the time-interval.
Counterexamples generated by waves concentrated on the support of the source may not arise in 1-d since  light rays hitting the boundary are only of hyperbolic type. 

Going back to the multi-d case under consideration, the lack of observability inequalities of the form \eqref{obs-loss}  shows that, necessarily,  an infinite number of derivatives may be lost on the measurement of the sources $g$, and thus, one has to impose some added restrictions on them to prevent concentration phenomena like \eqref{invisible} (see the pseudo-differential condition in assumption A3 below).

Within this class of sources $g$, the sidewise observability inequality \eqref{obs} will be proved under a  microlocal geometrical condition (see assumption A2 below), inspired (but different !) from the Geometric Control Condition introduced in \cite{B-L-R}. Roughly, it guarantees that all rays emanating from the support of the source reach the observation region without earlier bouncing on the support of the source.
This condition is sharp in terms of the geometry of the support of the sources $O$  and the measurement subset  $O'$ and also in what concerns the sidewise observability time. 

\subsection{Extensions and open problems.}
The methods of this paper could be employed to handle other related problems such as:
\begin{itemize}
\item The simultaneous initial and boundary source sidewise observation. We refer to \cite{Zuazua} for a complete analysis in 1-d.
\item The problem treated in \cite{Baudoin} where, on an annular domain $\Omega = A(R_{1},R_{2}) = \{x \in \R^{n}, \, R_{1} < \vert x \vert < R_{2}\}$ of $\R^{n}$, initial data are observed out of measurements on the exterior part of the boundary, under suitable conditions on the sources with support on the interior  boundary.
\end{itemize}

Similar questions on the sidewise boundary observability and source identification are also of interest for other models such as, for instance Schr\"odinger, plate and heat equations, the elasticity system and thermoelasticiy, all of them rather well understood in the control of classical boundary control. But their analysis would require of significant further developments.

\subsection{Structure of the paper}
The paper is organized as follows. In Section \ref{sec.2} we state the main results, and 
Section \ref{sec.3} is devoted to present some preliminary results. Most of the tools presented here are classical and we recall them in order to standardize the notations and  make the paper self-contained. We start with  the geometrical setting  and we present in particular the generalized bicharacteristic curves and the partition of the cotangent space of the boundary $T^{*}\d\L$. We also introduce the spaces of pseudo-differential symbols that will play the role of test functions on which we build the microlocal defect measures, of great importance in the proof. In Section \ref{sec.4}, we present a geometric consequence of Assumption A2 and we perform a pseudo-differential multiplier calculus up to the boundary, in the spirit of \cite{Lions}, that will play a central role in  the proof. 
Section \ref{sec 5} is  mostly devoted to the proof of the main result namely Theorem \ref{theo}. In Section \ref{sec 6}, we  present the proof of Theorem \ref{Obs-Loss-1}, essentially based on the microlocal behavior of the solutions to \eqref{waveequation}.  We also  present the proof of  Corollary \ref{Obs-Loss-3} where we construct a sequence of boundary data $(g_{k})$  concentrating on the glancing set.

\noindent\textbf{Acknowledgements.}
The authors thank Nicolas Burq for fruitful discussions about concentration of waves  near elliptic and glancing points of the boundary. The authors also thank Nicola de Nitti for his help on designing and executing the figures of the paper.

The research of the first author was partially supported by the Tunisian Ministry for Higher Education and Scientific Research  within the LR-99-ES20 program.
The second author has been funded by the Alexander von Humboldt-Professorship program, the Transregio 154 Project ``Mathematical Modelling, Simulation and Optimization Using the Example of Gas Networks" of the DFG, the ModConFlex Marie Curie Action, HORIZON-MSCA-2021-$d$N-01, the COST Action MAT-DYN-NET, grants PID2020-112617GB-C22 and TED2021-131390B-I00 of MINECO (Spain), and by the Madrid Government -- UAM Agreement for the Excellence of the University Research Staff in the context of the V PRICIT (Regional Programme of Research and Technological Innovation).

\section{Statement of the results}\label{sec.2}
\subsection{Sidewise observability}\label{sec.2.1}

Let $\Omega$ be a domain of $\R^{n}$ admissible in the sense of assumption A1, and  $O$ a subset of the boundary $\d\Omega$ strictly concave. And consider $O'$ a subset of $\d\Omega$ such that $\overline{O}\cap \overline{O'} = \emptyset$. We start with the geometric condition we will impose to the pair $\{O,O'\}$.

First, we recall  that given the cylinder $\L = \R\times \Omega$ with $\Omega$ of class $\mathcal{C}^{\infty}$,  we can define the Melrose-Sj\"ostrand   compressed cotangent bundle of $\L$,  $T^{\ast}_{b}\L= T^{\ast}\L \cup T^{\ast}\d\L$. In addition, the matrix $A=(a_{ij}(x))$  being also of class $\mathcal{C}^{\infty}$, we have a flow on $T^{\ast}_{b}\L$, constituted of generalized bicharacteristic curves of the wave operator
\begin{equation*}
P_{A}=\partial _{t}^{2}-\sum_{i,j=1}^{n}\partial_{x_{j}}(a _{ij}(x)\partial_{x_{i}}) ,
\end{equation*}
 the celebrated Melrose-Sj\"ostrand flow (see \cite{MeSj}). We refer the reader to Section \ref{boundary-geometry} for further details and precise definitions of these facts.

In particular, we recall the partition of the cotangent bundle of the boundary $T^{\ast}\d\L$ into elliptic, hyperbolic and glancing sets :
\begin{equation}
T^{\ast}\d\L = \E \cup \H \cup \G .
\end{equation}

Now, consider  an open subset  $\mathcal{O}$ of $\d\Omega$, strictly concave in the sense of assumption A1, such that $\overline{O}\subset \mathcal{O}$ and $\overline{\mathcal{O}}\cap \overline{O'}= \emptyset$. One can easily check that this is possible since A1  is an open condition. 

\smallskip
\textbf{Assumption A2: SGCC }
\newline
We assume that there exists a time $T_{0}>0$  such that  every generalized bicharacteristic curve issued from the boundary $\mathcal{O}$ at $t=0$, intersects the boundary $O'$ at a strictly gliding point , without intersecting  $\overline{\Gamma}$, and before the time $T_{0}$.
 \begin{figure}[h!]\label{bichar}
 \begin{center}
 \subfigure
{\includegraphics[scale=0.9]{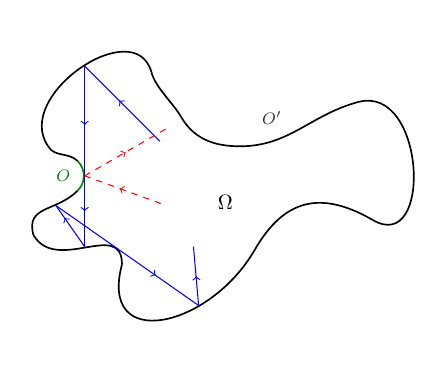}}

 \caption{Bicharacteristic rays passing throw $O$}
\end{center}
\end{figure}

 \begin{rema}\label{remA2} 
 \begin{enumerate}
 
 \item The definition of  strictly gliding point of the boundary will be given in Section \ref{boundary-geometry}.
 
 \item The notation (SGCC) stands for \textit{sidewise geometric control condition}. In what follows, we provide some precisions. 
 
 \item  Set $\mathcal{U}=\R\times \mathcal{O}$. The generalized bicharacteristic curves issued from points of the boundary $\mathcal{U}$ are of two types and can be described through their projection on the basis, i.e the $(t,x)-$space. On one hand we have the curves that are transverse to $\d\L$ and in this case we have  two hyperbolic fibers issued from the same hyperbolic point $m_{0}\in \d\L$. At $m_{0}$, we have a hyperbolic reflection. On the other hand, the curve is tangent to $\d\L$ at $m_{0}$ ( one order tangency ) and lies in $\L =\R\times \Omega$ otherwise. In the latter case, the generalized bicharacteristic curve can be interpreted as a ``free bicharacteristic curve'' since it's an integral curve of the hamiltonian field attached to the wave symbol ( see Section \ref{boundary-geometry}). 
 
Condition (SGCC) requires that each one of these curves  starting from $\mathcal{U}$ at $t=0$, to intersect the boundary $\Gamma'$ at a strictly gliding point , without intersecting  $\overline{\Gamma}$, and before the time $T_{0}$. In this sense, this condition is stronger than the classical (GCC) of Bardos, Lebeau and Rauch \cite{B-L-R} that needs  the rays to hit $\d\Omega$ at non diffractive points.

 \item For instance if  $\gamma = \gamma(s)$ is a ray issued from $\mathcal{U}$, we  have 
 $\gamma(0) = \rho \in T_{b}^{\ast}\L_{\vert\mathcal{U}}$, $\gamma(s_{0}) =\rho_{1}\in T_{b}^{\ast}\L_{\vert\Gamma'}$ for some $s_{0} \in ]0,T_{0}[$, where $\rho_{1}$ is a strictly gliding point, and moreover $\gamma(s) \notin  T_{b}^{\ast}\L_{\vert\overline{\Gamma}}$
for $0<s< s_{0}$.
\newline
In particular we can allow $\gamma(s)$ to live on the boundary, outside $T_{b}^{\ast}\L_{\vert\overline{\Gamma}}$ for some values of $s \in  ]0,s_{0}[$.

\item Notice that we don't make any assumption on the rays  that don't intersect the open set $\mathcal{U}$ of the boundary. From this point of view, (SGCC) is weaker than the classical condition (GCC). 
\item  Remark that if $O$ is strictly convex, then obviously, (SGCC) cannot be satisfied ( see Fig.\ref{convex}). Therefore, assumption A1 seems to be a well adapted framework to set up the microlocal condition A2.

\begin{figure}[h!]
 \begin{center}\label{convex}
 \includegraphics[scale=0.7]{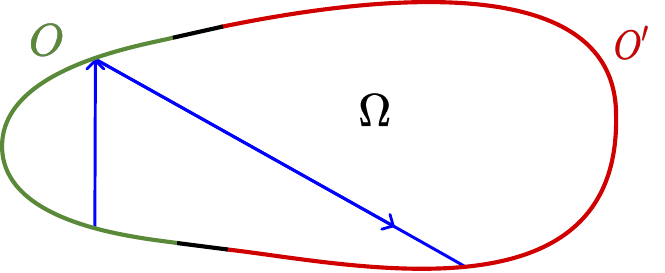}
\caption{Convex boundary. In blue, a geodesic ray.}
 \end{center}
 \end{figure}

 \end{enumerate}
\end{rema}

Finally, we introduce the last assumption, namely a boundary condition on the data $g$. For this purpose, we recall that the lateral boundary $\d\L$ of the cylinder $\L=\R\times \Omega$ is a  submanifold of $\R^{n+1}$, of dimension $n$ and class $\mathcal{C}^{\infty}$. We will denote by $(t,x')=(t,x'_{1},...,x'_{n-1})$ a system of local coordinates on $\d\L$.

\smallskip
\textbf{Assumption A3: Boundary condition fulfilled by observable  sources}

We assume  one of the  following conditions :

\textbf{A3.a} There exists a  polyhomogeneous pseudo-differential operator $B_{\alpha}= b_{\alpha}(t,x';D_{t}, D_{x'})$  on $\d\L$, of order $\alpha >0$, such that $CharB_{\alpha} \subset \H$ and
\begin{equation} \label{regularity-0}
b_{\alpha}(t,x';D_{t}, D_{x'})g =0.
\end{equation}

\textbf{A3.b} There exists a family of polyhomogeneous pseudo-differential operators $c_{\alpha}(t,x';D_{x'})$ in the $x'$-variable on $\d\L$, smooth with respect to respect to $t$, elliptic  of order $\alpha >0$ such that
\begin{equation} \label{regularity}
c_{\alpha}(t,x';D_{x'})g  = 0.
\end{equation}

\textbf{A3.c} There exists $\mathcal{U}_{M}$ an open neighborhood of $\overline{\Gamma}_{M}$ in $\d\L$,  there exists $\alpha >0$ and a constant $C_{\alpha}>0$ such that for every $u$ solution of system \eqref{waveequation},  the boundary trace 
\newline
$(\d_{n}u+\d_{t}u)_{\vert \d\L}$ satisfies 
\begin{equation} \label{trace regularity}
\Vert (\d_{n}u+\d_{t}u)_{\vert \d\L}\Vert_{H^{\alpha}(\mathcal{U}_{M})}  \leq C_{\alpha}\Vert g\Vert_{H^{1}(\Gamma_{M})} .
\end{equation}

\begin{rema}
For the definition of polyhomogeneous pseudo-differential operators on $\d\L$, see Section \ref{pdo}. In particular, we recall that the characteristic set of  $B_{\alpha}= b_{\alpha}(t,x';D_{t}, D_{x'})$  is given by 
$$
CharB_{\alpha}=\{(t,x';\tau, \xi')\in T^{\ast}\d\L , \,\, \sigma(b_{\alpha})(t,x';\tau, \xi') =0\}
$$
where $\sigma(b_{\alpha})$ is the principal symbol of $B_{\alpha}$.
\end{rema}
We are now ready to state our main theorem.

\begin{theorem} \label{theo}
Under assumptions A1, A2 and A3,  for every $T>T_{0}$, there exists  $C > 0$ such that every solution of  \eqref{waveequation}, satisfies the observability estimate  
\begin{equation}\label{obs-theo}
\Vert g\Vert_{H^{1}(\Gamma_{M})} \leq C \Vert \d_{n}u_{\vert\Gamma'}\Vert_{L^{2}(\Gamma'_{M+T} )}.
\end{equation}
\end{theorem}

\begin{rema}\label{remark}
\begin{enumerate}
\item In case assumption A3.a is satisfied, we can relax assumptions  A1 and A 2. Indeed, we may only assume  the subset $O$ of the boundary $\d\Omega$ to be concave and not necessarily strictly concave. In particular, it can be locally a hyperplane. In addition, we may assume  A2  only for transverse ( hyperbolic ) rays.

\item Condition A3.b ensures some \`a priori spatial regularity on the data $g$,  yielding  microlocal  regularity of $g$ near the elliptic and the glancing sets of the boundary. For instance, it is fulfilled if  $g$ doesn't depend on the space variable $x'$, i.e $g=g(t)$. In the same spirit, if we assume 
\begin{equation*}
\Vert\nabla_{x'}u_{\vert \d\L}\Vert_{H^{\alpha}(\mathcal{U}_{M})}  \leq C_{\alpha}\Vert g\Vert_{H^{1}(\Gamma_{M})} ,
\end{equation*}
for some $\alpha >0$, we get the same positive conclusion, as a byproduct of the previous argument .

\item In Assumption A3.c , the open set $\mathcal{U}_{M}$ can be taken in the form $(-\varepsilon,M+\varepsilon)\times \mathcal{O}$, where $\mathcal{O}$ is an open neighborhood of $\overline{O}$ in $\d\Omega$. This condition can be interpreted as a conditional stability assumption. See for instance V. Isakov \cite{Isakov}.

\item Obviously, the three conditions  a), b) and c)  of Assumption A3 are each of them sufficient and complementary. One could consider other assumptions guaranteeing the conclusion of Theorem \ref{theo}.

\item In the setting of assumption A3.a, one can  for instance, consider the case where the boundary data $g$ is subject to a wave equation. With $\chi = \chi (t,x) \in \mathcal{C}^{\infty}_{0}(\Gamma_{M})$, consider the system
\begin{equation}
\left\{ 
\begin{array}{c}
P_{A}u=\partial _{t}^{2}u-\sum_{i,j=1}^{n}\partial_{x_{j}}a _{ij}(x)\partial_{x_{i}}u=0\quad \text{in } \L 
\\ 
\\
u(t,.)= \chi (t,x) g(t,.)\quad \text{on } \d\L 
\\
\\
P'_{A}g = \partial _{t}^{2}g- \beta \sum_{i,j=1}^{n-1}\partial_{x'_{j}}a _{ij}(x',0)\partial_{x'_{i}}g=0\quad \text{on } \d\L
\\
\\
u(0,.)=\partial _{t}u(0,.)= 0 \quad \text{on  }  \Omega
\\
\\
g(0,.)=g_{0}\in H^{1}(\d\L), \quad \text{and}\quad\partial _{t}g(0,.)=g_{1}\in L^{2}(\d\L)
\end{array}%
\right.  \label{waveequation-1}
\end{equation}
where $\beta >0$. One can easily check that assumption A3.a is fullfilled as soon as $\beta >1$.

However, if $\beta \leq 1$, the  characteristic set of $P'_{A}$ is contained in the union $\E \cup \G$ of the elliptic set  and the glancing set. In this case, one can construct  a sequence of sources $(g_{k})$ such that the corresponding sequence of solutions $(u_{k})$ to system    \eqref{waveequation-1} violates the observability estimate \eqref{obs-theo}, with a loss of compactness  located in  $\E$ or $\G$, see Theorems \ref{Obs-Loss-1}  and \ref{Obs-Loss-3}.
\item To summarize:  Even if, thanks to (SGCC), we can microlocally control the source $g$ near the hyperbolic set of $\d\L$,   it still may develop singularities on the elliptic set, and/or travelling along some characteristic curves of the glancing set. In fact, as we will see in the proof of Theorem \ref{theo}  the analysis on these sets requires a special attention. Assumption A3.a , A3.b or A3.c  above are set to  insure additional regularity on $g$ that avoids the  rising of such singularities. 
\end{enumerate}
\end{rema}

\subsection{On the lack of sidewise observability}

We present now the results concerning the lack of observability, even in the weaker version \eqref{obs-loss}. These  negative results  ensure a loss of an infinite number of derivatives for all possible geometric configurations. Here we do not need any of the geometric conditions A1 or A2, that is, we work on a general bounded and smooth domain $\Omega$ and any partition of its boundary.

The proofs of these results will be given in  Section \ref{sec 6}.

\begin{theorem} \label{Obs-Loss-1}
For every $s < 0$, there exists a sequence of sources $(g_{k})_{k\geq 1} \subset H^{1}(\d\L)$ supported in $\overline{\Gamma}_{M}$, such that the solutions $(u_{k})$ of system \eqref{waveequation} satisfy 
\begin{equation}\label{obs-Loss-2}
\mathop{\lim}_{k \to \infty}\Vert g_{k}\Vert_{H^{s}(\Gamma_{M})}=1 \quad \text {and} \quad  \mathop{\lim}_{k \to \infty}\Vert \partial_{n}u_{{k}_{\vert\d\Omega}} \Vert_{L^{2}(\Gamma'_{M+T})} = 0 ,
\end{equation}
for every $T>0$.
In particular, the lack of compactness of the sequence $(g_{k})$ in $H^{s}(\Gamma_{M})$ is located in the elliptic set $\E$ of the boundary.
\end{theorem} 

And we deduce from this theorem :

\begin{corollary} \label{Obs-Loss-3}
For every $s < 0$, there exists a sequence of sources $(g_{k})_{k\geq 1} \subset H^{1}(\d\L)$ supported in $\overline{\Gamma}_{M}$, such that the solutions $(u_{k})$ of system \eqref{waveequation} satisfy \eqref{obs-Loss-2}
for every $T>0$.

In particular, the lack of compactness of the sequence $(g_{k})$ in $H^{s}(\Gamma_{M})$ is located in the glancing set $\G$ of the boundary.

\end{corollary}

\begin{rema} 
Actually, as we will see in the proof (cf. Section \ref{sec 6}), we choose a sequence $(g_{k})$ supported in $\overline{\Gamma}_{M}= [0,M]\times \overline{O}$ such that  for some fixed $\alpha >1$,  $\Vert g_{k}\Vert_{H^{\alpha}}$ is bounded outside the elliptic set $\E$ of the boundary. The propagation of the $H^{\alpha}$-wave front will then provide the desired result. 
\end{rema} 

\begin{rema} 
In view of Theorem \ref{Obs-Loss-1},  we can not expect the sidewise observability estimate \eqref{obs-theo} to hold, unless an infinite number of derivatives is lost. Therefore, in order to get sidewise observability estimates in Sobolev norms,  structural conditions on the sources need to be imposed,  such as those of assumption A3.
%
%

\end{rema}

\begin{rema} 
To close this section and before going into the proofs, let us summarize the strategy one should follow to obtain sidewise observability for system \eqref{waveequation}.

First, we have to adress the problem only on well designed domains $\Omega$, i.e those satisfying assumption A1. Secondly, we  choose the measurements domain, i.e a subset $O'$ of the boundary $\d\Omega$, $\overline{O} \cap \overline{O'}=\emptyset$,  as sharp as possible, such that (SGCC) is fullfilled. For instance, in the case of the annular domain ( Fig.\ref{geometryO}), if $O$ is the interior boundary, then $O'$ is the exterior boundary. And finally, we make sure that the  boundary source $g$ we aim to observe is admissible, i.e it satisfies some \`a priori condition in the spirit of condition A3, that prevents the presence of invisible solutions. 
\end{rema}

\section{Some Geometric Facts, Operators and Measures}\label{sec.3}
\subsection{Geometry }\label{geo}
Near a point  $m_{0}$ of the boundary $\d\Omega$, taking advantage of the regularity of $\Omega$, we can define a system of geodesic local coordinates $x=(x_{1},x_{2},....,x_{n}) \longrightarrow y=(y_{1},y_{2},....,y_{n})$ such that 
$$
\Omega= \{ (y_{1},y_{2},....,y_{n}), \,\, y_{n}>0\}, \quad \d\Omega = \{ (y_{1},y_{2},....,y_{n-1},0)\}=\{ (y',0)\}
$$
where the wave operator is given by
$$
P_{A} = -\partial_{t}^{2}+\Big(\partial_{y_{n}}^{2}+\sum_{1\leq i,j \leq n-1}\partial_{y_{j}}b_{ij}(y)\partial_{y_{i}}\Big) +M_{0}(y)\partial_{y_{n}} + M_{1}(y,\partial_{y'}) .
$$
Here, the matrix $(b_{ij}(y))_{ij}$ is of class $\mathcal{C}^{\infty}$, symmetric,  uniformly definite positive on a neighborhood of $m_{0}$, $M_{0}(y)$ is a real valued function of class $\mathcal{C}^{\infty}$, and $ M_{1}(y,\partial_{y'})$ is a tangential differential operator  of order $1$ with $\mathcal{C}^{\infty}$ coefficients.

In the sequel, we will come back to the notation $(t,x)=(t,x',x_{n})=(t,y',y_{n})$, and we shall write
$$
P_{A} = \partial_{n}^{2}+R(x_{n},x',D_{x',t}) + M_{0}(x)\partial_{n} + M_{1}(x,\partial_{x'})
$$
Notice that, in this coordinates system,  the principal symbol of the wave operator $P_{A}$ is given by 
$$
\sigma(P_{A}) = -\xi_{n}^{2}+r(x,\tau,\xi') = -\xi_{n}^{2} +\Big(\tau^{2}-\sum_{1\leq i,j \leq n-1}a_{ij}(x)\xi_{i}\xi_{j}\Big).
$$
We shall set $r_{0}(x',\tau,\xi')=r(x',0, \tau,\xi')$ and we denote $m_{1}=m_{1}(x,\xi')$ the symbol of the vector field $M_{1}$.

\subsection{Generalized bicharacteristic rays}\label{boundary-geometry}

Let us introduce the compressed cotangent bundle of Melrose-Sj\"ostrand 
$T^*_b\L=T^{\ast }\L\cup T^{\ast }\partial \L$. We recall that we have a natural projection
\begin{equation}
\label{projection1}
\pi\ :\ T^{*}\R^{n+1}\mid _{\overline \Omega} \,\rightarrow T^{*}_{b}\L
\end{equation}

and we equip $T^{\ast }_b\L$ with the induced topology. 

 Given the matrix $A(x)=(a_{ij}(x))$, we denote by $p_A(x;\tau,\xi)= \tau^2-^t\xi A(x)\xi$, the principal symbol of the wave operator, and  
$$
\text{Char}(P_{A}) =  \{(t,x;\tau,\xi), p_{A}(x,\tau, \xi) = \tau^2-^{t}\xi A(x)\xi = 0\},
$$ 
the characteristic set, and  $\Sigma_{A} = \pi (Char(P_{A}))$.
In addition, we recall  the hamiltonian field associated to $p_A$
$$
H_{p_A} = 2\tau \partial_t -2 ^t\xi A(x)\partial_x +\sum_{k=1}^{n} \, ^t \xi\partial_{x_k}A(x)\xi\partial_{\xi_k} .
$$

Also, we recall the following partition of $T^{\ast}(\d\L)$ into elliptic, hyperbolic and glancing sets:
 
 \begin{equation}\label{Boundary}
 \#\Big\{\pi^{-1}(\rho) \cap Char(P_A)\Big\} = 
\left\{ 
\begin{array}{c}
0 \quad if  \quad \rho \in \E
 \\ 
1 \quad if  \quad \rho \in \G
\\ 
2 \quad if  \quad \rho \in \H
\end{array}%
\right.  
\end{equation}
For the sake of simplicity, we will develop the rest of this section in a system of local geodesic coordinates as introduced in section \ref{geo}.  We recall that we have locally

$$
\L = \{(t,x) \in \R^{n+1}, \,x_{n}>0\}  \quad \text{and} \quad \d \L = \{(t,x) \in \R^{n+1}, \, x_{n}=0\} .
$$

We also get :
$$
\E = \{r_{0}<0\}, \quad\quad\H = \{r_{0}> 0\}, \quad\quad\G = \{r_{0}=0\} .
$$

Notice that using the projection $\pi$, one can identify the glancing set $\G$ with a subset of $T^{\ast }\R^{n+1}$.

\begin{defi}\label{def:nondiffractive}
\begin{enumerate}
\item A point $\rho \in T^{\ast }\partial \L\backslash 0$ is 
nondiffractive  if $\rho \in \H$ or if $\rho \in \G$ and the free bicharacteristic $%
(\exp sH_{p_A})\widetilde{\rho }$ passes over the complement of $\overline{\L}$ for
arbitrarily small values of $s,$ where $\widetilde{\rho }$ is the unique
point in $\pi^{-1}(\rho)\cap Char(P_A)$. 

\item $\rho \in T^{\ast }\partial \L\backslash 0$ is strictly gliding if $\rho \in \H$ or if $\rho \in \G$ and $H_{p_A}^2(x_{n})(\rho) < 0$.

In the latter case, the  free \textit{bicharacteristic }ray $\gamma$ issued from $\rho$ leaves the boundary $\partial \L$ and enters in $T^{\ast }(\R^{n+1}\setminus \overline{\L})$ at  $\widetilde{\rho }=\pi^{-1}(\rho)$.

\item $\rho \in T^{\ast }\partial \L\backslash 0$ is strictly diffractive  if $\rho \in \G$ and $H_{p_A}^2(x_{n})(\rho) > 0$.

This means that there exists $\varepsilon >0$ such that $(\exp sH_{p_A})\widetilde{\rho } \in T^{\ast }\L$ for $0<\vert s\vert <\varepsilon$.

\end{enumerate}
\end{defi}
\begin{defi}\label{def:diffractive} 

We shall denote by $\G_d$ the set of strictly diffractive points and by $\G_{sg}$ the set of strictly gliding points.
\end{defi}

\begin{rema}\label{rem:diffractive}
\begin{enumerate}

\item  Under assumption A1, we notice that over $\Gamma$,  the glancing set $\G$ is reduced to $\G_{d}$,  i.e
$$
\G_{\vert \Gamma} \subset \G_{d}.
$$
Namely all generalized bicharacteristic curves issued from points of $\G_{\vert \Gamma}$  have a first order tangency with the boundary .

\item In local geodesic coordinates, the sets $\G_{d}$ and $\G_{sg}\setminus \H$ are given by 
\begin{equation}\label{diffractive}
\G_{d}= \{\xi_{n}= r_{0}=0, \, \partial_{n}r_{\vert x_{n}=0}>0\}, \quad\quad \text{and} \quad\quad \G_{sg}\setminus \H= \{\xi_{n}= r_{0}=0, \, \partial_{n}r_{\vert x_{n}=0} < 0\} .
\end{equation}

\end{enumerate}
\end{rema}

\begin{defi}\label{def:gene-bichar}
A generalized bicharacteristic ray is a continuous map
$$
\R \supset I \setminus B\ni s \mapsto \gamma (s) \in T^{\ast }\L \cup \G \subset T^{\ast }\R^{n+1}
$$
where $I$ is an interval of $\R$, B is a set of isolated points, for every $s  \in I \setminus B$, $\gamma (s) \in \Sigma_A$ and $\gamma$ is differentiable as a map with values in $T^{\ast }\R^{n+1}$, and 
\begin{enumerate}
\item If $\gamma (s_0) \in T^{\ast }\L \cup \G_d$ then $\dot{\gamma}(s_{0})=H_{p_A}(\gamma(s_{0}))$.

\item If $\gamma (s_0) \in \G \setminus \G_d$ then $\dot{\gamma}(s_0)=H_{p_A}^G(\gamma(s_0))$, where  $H_{p_A}^G= H_{p_A} + (H_{p_A}^2x_{n} /H_{x_{n}}^2p_A)H_{x_{n}}$.

\item For every $s_0 \in B$, the two limits $\gamma(s_0 \pm 0)$ exist and are the two different points of the same hyperbolic fiber of the projection $\pi$.
\end{enumerate}
\end{defi}

\begin{rema}
\begin{enumerate}

\item We recall that if  $\Omega$ has no contact of infinite order  with its tangents, the Melrose-Sj\"ostrand flow is globally well defined. 

\item In the interior, i.e in $T^{*}\L$, a generalized bicharacteristic  is simply a classical bicharacteristic ray of the wave operator whose projection on the basis is a geodesic of $\Omega$ equipped with the metric $(a^{ij})=(a_{ij})^{-1}$ .

\item  Finally, $\gamma$ can be considered as a continuous map on the interval $I$ with values in $T^{\ast }_b\L$.

\end{enumerate}
\end{rema}

\subsection{Pseudo-differential operators }\label{pdo}
In this section, we introduce the classes of pseudo-differential operators we shall use in this paper. We start with the operators on the cylinder $\L$.

Let $\mathcal{A}$ be the set of pseudo-differential operators of the form $Q=Q_{i}+Q_{\d}$ where $Q_{i}$ is a classical pseudo-differential operator ,  compactly supported in $\L$ and $Q_{\d}$ is a classical  tangential pseudo-differential operator, compactly supported near $\d\L$. More precisely, $Q_{i}=\varphi Q_{i}\varphi$ for some $\varphi \in \mathcal{C}^{\infty}_{0}(\L)$ and $Q_{\d}=\psi Q_{\d}\psi$ for some $\psi(t,x_{n}) \in \mathcal{C}^{\infty}(\R \times ]-\alpha, \alpha[)$. $\mathcal{A}^{s}$ will denote the elements of $\mathcal{A}$ of order s.

On the other hand, the boundary $\d\L = \R\times \d\Omega$ is a smooth manifold of dimension $n$ without boundary. Following L.H\"ormander \cite{Hormander-III} and using a system of local charts,   we can define for  $m \in \R$, the space of polyhomogeneous pseudo-differential operators  $\Psi_{phg}^{m}(\d\L)$ on $\d\L$, associated with symbols in $S_{phg}^{m}(T^*\d\L)$. These operators enjoy all classical properties of continuity and composition.

\subsection{Microlocal defect measures }
\label{mdm}
Here we use notations of section \ref{boundary-geometry}. Denote
 \begin{equation*}
\left\{ 
\begin{array}{c}
Z=\pi(CharP_{A}) , \quad\quad \hat{Z} = Z \cup \pi( T^{\ast}\overline{\L}_{\vert x_{n}=0}) ,
\\
 \\ 
SZ = (Z\setminus \overline{\L})/ \R_{+}^{\ast} , \quad\quad S\hat{Z} = (\hat{Z}\setminus \overline{\L})/ \R_{+}^{\ast} .
\end{array}%
\right.  
\end{equation*}
and for $Q \in \mathcal{A}^{0}$ with principal symbol $\sigma(Q)=q$, set 
$$
\kappa(q)(\rho)=q(\pi^{-1}(\rho)).
$$
We define also for $u \in H^{1}(\L)$
$$
\phi(Q,u) = (Qu,u)_{H^{1}} = \int_{\L}\Big(\nabla_{t,x}Qu .\nabla_{t,x}\overline{u} + Qu.\overline{u} \Big)dxdt .
$$
Finally, let $(u_{k})$ be a sequence of functions weakly converging to $0$ in $H^{1}_{loc}(\L)$. In \cite{Lebeau} and \cite{BurqLebeau}, the authors  prove the following result:
\begin{theorem}[Burq-Lebeau \cite{BurqLebeau}]
\label{measuresB-L}
There exists a subsequence of $(u_{k})$ (still denoted by $(u_{k})$)  and a positive Radon measure $\mu$ on $S\hat{Z}$ such that
\begin{equation*}
\lim_{k\rightarrow \infty}\phi(Q, u_{k}) = \langle  \mu , \kappa (q)  \rangle , \quad\quad  \forall Q \in \mathcal{A}^{0} . 
\end{equation*}
\end{theorem}
We will refer to $\mu$ as a microlocal defect measure associated to the sequence $(u_{k})$.

On the other hand, on the boundary $\d\L$, we can make use of the classical notion of microlocal  defect measure introduced by P. G\'erard in \cite{Ge91}. More precisely, for every sequence of functions $(v_{k})$   weakly converging to $0$ in $H^{1}_{loc}(\d\L)$, there exists a positive Radon measure $\tilde{\mu}$ on $S^{\ast}(\d\L)$ such that we have,  up to a subsequence
\begin{equation*}
\lim_{k\rightarrow \infty}(Qv_{k}, v_{k})_{L^2(\d\L)} = \langle  \tilde{\mu} , \vert\eta\vert^{-2}\sigma (Q)  \rangle , \quad\quad  \forall Q \in \Psi_{phg}^{2}(\d\L) . 
\end{equation*}
Here we have denoted by $(y,\eta)$ the standard element of $T^{\ast}(\d\L)\setminus0$.

We will remind the properties of these measures in some steps of the proof later, see Section \ref{measures-prop}.

\section{Preliminary results}\label{sec.4}

\subsection{A Geometric Lemma}\label{subsec.4}
\begin{figure}[h!]\label{half-bichar}
 \subfigure
{\includegraphics[scale=0.85]{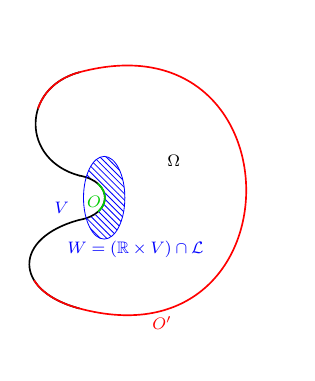}}
 \quad\quad\quad\quad
  \subfigure
  {\includegraphics[scale=0.3]{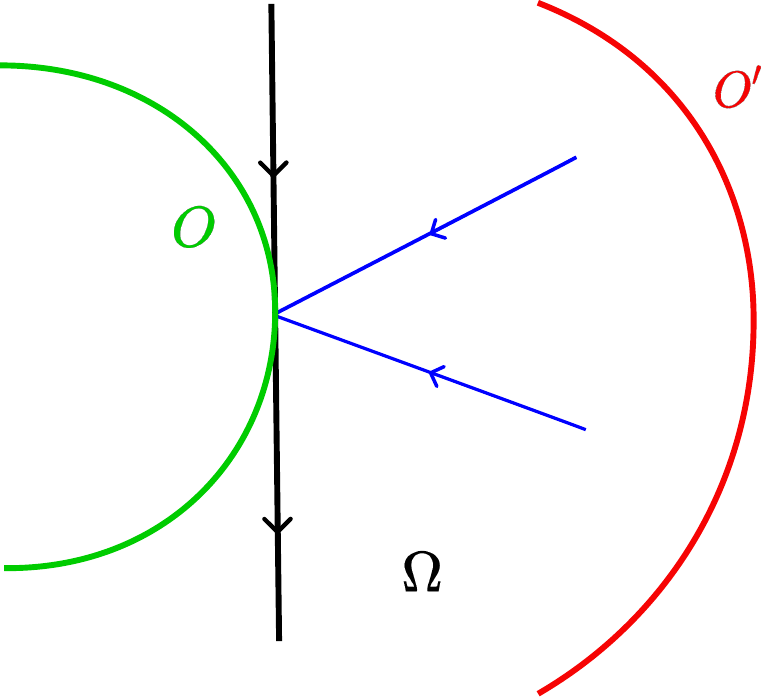}}
  \caption {On the left, interior neighborhood of $\Gamma$. 
  \\
  On the right, tangent ( black ) and hyperbolic ( blue ) half bicharacteristic rays}
\end{figure}

Let $O$ (resp. $\mathcal{O}$) be the open subset of $\d\Omega$ introduced in the statement of Assumption A1 (resp. A2 ),  and set $\mathcal{U} =\R\times \mathcal{O}$. Consider $V$ a neighborhood of $\overline{O}$ in $\R^{n}$ such that $V \cap \d\Omega \subset \mathcal{O}$. $\R \times V$ is an open neighborhood of $\overline{\Gamma} = \R\times \overline{O}$ in $\R^{n+1}$. In this setting $W = \R\times (V\cap \Omega)=(\R\times V)\cap \L$ is an interior neighborhood of the boundary $\overline{\Gamma}$ ( see Figure \ref{half-bichar}).  
On the other hand, consider $\rho \in T^{*}W \cap Char(P_{A})$  and denote $\gamma = \gamma(s)$ the generalized bicharacteristic issued from $\rho$, i.e $\gamma(0)=\rho$. In addition, we define by   $\gamma^{+}=\{\gamma(s), s>0\}$ , resp. $\gamma^{-}=\{\gamma(s), s<0\}$   the outcoming half bicharacteristic and the incoming half bicharacteristic at $\rho$, see Figure \ref{half-bichar}.

\begin{lemma}\label{bound-nbhood}
 With the notations above and under assumptions A1 and  A2, for every $T>T_{0}$, there exists  $V$ neighborhood  of $\overline{O}$ in $\R^{n}$, $V\cap \d\Omega \subset \mathcal{O}$,  such that for every $\rho \in T^{*}(W)\cap Char(P_{A})$, one of the two half bicharacteristics issued from $\rho$, the outcoming one or the incoming one, travelling at speed one, intersects the boundary $\Gamma'$ at a strictly gliding point , without intersecting the boundary $\overline{\Gamma}$, and before the time $T$.
 
 We will say that this half bicharacteristic satisfies (SGCC).
\end{lemma} 

\begin{proof} 
For $\rho \in T^{*}W\cap Char(P_{A})$, denote by $\gamma_{\rho}=\{\gamma_{\rho}(s), s \in \R\}$ the generalized  bicharacteristic issued from $\rho$. In particular, $\gamma_{\rho}(0) = \rho$. Assume that $\gamma_{\rho}$ intersects $\mathcal{U}$ for some value $s_{1}<0$ at a hyperbolic or at a glancing point.  According to assumption A2, we then get that for some $s \in \R$ such that $s-s_{1} < T_{0}$,   $\gamma_{\rho}(s)$ is a strictly gliding point of the boundary $\Gamma'$ and, in addition $\{\gamma_{\rho}(s'), s_{1} <s' <s\} \cap T^{*}_{b}\L_{\vert \overline{\Gamma}} = \emptyset$. In this case, we see that the statement of Lemma \ref{bound-nbhood} is satisfied by the outcoming half bicharacteristic issued from $\rho$. Obviously, the case $s_{1}>0$ can be treated in a similar way. According to this, we may only focus on the points $\rho$ close to $\overline{\Gamma}$ such that $\gamma_{\rho}=\{\gamma_{\rho}(s), s \in \R\}$ doesn't intersect $\overline{\Gamma}$ for $s \in ]-T_{0},T_{0}[$. In addition, due to the compactness of $\overline{O}$, it suffices to prove that  every glancing point $\rho \in \G_{\vert\mathcal{U}} \subset T^{*}\d\L_{\vert\mathcal{U}}$  admits a neighborhood $V_{\rho}$ in $T^{*}(\R^{n+1})$ such that  conclusion of Lemma \ref{bound-nbhood} is valid for every $\rho' \in V_{\rho}\cap T^{*}\L$.

Before entering in the details of the proof, we warn the reader that if a generalized  bicharacteristic $\gamma_{\rho}$ hits the boundary transversally for some value $s_{0}$, that is at a hyperbolic point, we will denote this point by $\gamma_{\rho}(s_{0})$ , by abuse of notation. 

Consider then $\rho \in \G_{\vert\mathcal{U}} \subset T^{*}\d\L_{\vert\mathcal{U}}$  and let $s_{0}\in ]0,T_{0}[$  be a time such that the generalized  bicharacteristic $\gamma_{\rho}$ hits the boundary $\Gamma'$ at a strictly gliding point. Here we have two possibilities : a) $\gamma_{\rho}(s_{0})$ is a hyperbolic point or b) $\gamma_{\rho}(s_{0})$ a glancing strictly gliding point. We will discuss each one of these cases,  and in order to simplify the argument, we will work in local geodesic coordinates. 

\begin{itemize}

\item Case a) : $\gamma_{\rho}(s_{0})$ is a hyperbolic point. With the notations of Definition \ref {def:gene-bichar}, $s_{0} \in B_{\rho}$  where $B_{\rho}$ is a  set of isolated points in $\R$ such that the two limits $\gamma_{\rho}(s_0 \pm 0)$ exist and are the two different points of the same hyperbolic fiber of the projection $\pi$. Furthermore, we have 
\begin{equation}\label{hyp-0}
H_{p_{A}}x_{n}(\gamma_{\rho}(s_{0}-0))= \frac{dx_{n}}{ds}(\gamma_{\rho}(s_{0}-0)) = -2\xi_{n}(\gamma_{\rho}(s_{0}-0)) < 0 .
\end{equation}

Consequently, for $\varepsilon >0$ small enough, $\gamma_{\rho}(s_{0}-\varepsilon)$ is an interior point, moreover, the $x_{n}$ and $\xi_{n}$- coordinates  satisfy
\begin{equation}\label{hyp-1}
-2\xi_{n}(\gamma_{\rho}(s)) =  \frac{dx_{n}}{ds}(\gamma_{\rho}(s))\leq -c, \quad \forall s\in [s_{0}-\varepsilon ,s_{0}[, \quad\quad \text{for some} \quad c>0 .
\end{equation}
This yields 
\begin{equation}\label{hyp-2}
\xi_{n}(\gamma_{\rho}(s)) \geq c/2, \quad \forall s\in [s_{0}-\varepsilon ,s_{0}[ .
\end{equation}

In addition, we may assume that $0<x_{n}(\gamma_{\rho}(s_{0}-\varepsilon))<\eta$ for some $\eta >0$ to be chosen later. Now we fix $\varepsilon >0$. Taking  into account the continuity of the Melrose-Sj\"ostrand flow, it's clear that  for  $0< \alpha <\frac{1}{4}  x_{n}(\gamma_{\rho}(s_{0}-\varepsilon)) $, one can find $V_{\rho}$  a small enough  neighborhood of  $\rho$ in $T^{*}\R^{n+1}$, such that for all  $\rho' \in V_{\rho}\cap T^{*}\L \cap Char(P_{A})$, 

\begin{equation}\label{hyp-3}
\vert x_{n}(\gamma_{\rho}(s_{0}-\varepsilon)) - x_{n}(\gamma_{\rho'}(s_{0}-\varepsilon))\vert \leq \alpha ,
\end{equation}
and

\begin{equation}\label{hyp-4}
\xi_{n}(\gamma_{\rho'}(s)) \geq c' , \quad \forall s\in [s_{0}-\varepsilon ,s_{0}[,
\end{equation}
for some $c'>0$.  In particular, this means that $\gamma_{\rho'}(s_{0}-\varepsilon)$ is an interior point since

\begin{equation} 
x_{n}(\gamma_{\rho'}(s_{0}-\varepsilon)) \geq \frac{3}{4}x_{n}(\gamma_{\rho}(s_{0}-\varepsilon))>0 .
\end{equation} 

In addition, notice that estimate \eqref{hyp-4} is valid  as long as $x_{n}(\gamma_{\rho'}(s)) >0$, so possibly for $s \in ]s_{0}-\varepsilon, s_{0}+\beta[$, $ \beta >0$  small . Finally,

\begin{equation}\label{hyp-5}
\left\{ 
\begin{array}{c}
x_{n}(\gamma_{\rho'}(s)) \leq  x_{n}(\gamma_{\rho'}(s_{0}-\varepsilon)) -2c' (s-s_{0}+\varepsilon) 
\\
\\
\leq \frac{5}{4}x_{n}(\gamma_{\rho}(s_{0}-\varepsilon)) -2c' (s-s_{0}+\varepsilon) \leq \frac{5}{4}\eta -2c' (s-s_{0}+\varepsilon)
\end{array}%
\right. 
\end{equation}
Consequently, we obtain that $x_{n}(\gamma_{\rho'}(s))$ vanishes for some $s \geq s_{0}+\frac{5}{8c'}\eta-\varepsilon$, which means that the bicharacteristic ray $\gamma_{\rho'}$ leaves $\L$ at a hyperbolic point before the time  $T>T_{0}$,  as soon as  $\frac{5}{8c'}\eta-\varepsilon < T-T_{0}$ .

\item Case b) : $\gamma_{\rho}(s_{0})$ is a glancing strictly gliding point. According to Definition \ref{def:nondiffractive}, we know in this case that 
\begin{equation}\label{hyp-6}
x_{n}(\gamma_{\rho}(s_{0})) = r(\gamma_{\rho}(s_{0})) = 0 \quad\quad \text{and}\quad\quad \frac{\partial r}{\partial x_{n}}(\gamma_{\rho}(s_{0})) <0 .
\end{equation}
 Let then $B(\gamma_{\rho}(s_{0}), \varepsilon)$ be the open ball of $T^{*}\R^{n+1}$ with center $\gamma_{\rho}(s_{0})$ and radius $\varepsilon$. It's clear that for $\varepsilon$ and $c>0$ suitable, one has 
 
 \begin{equation}\label{hyp-7}
 \frac{\partial r}{\partial x_{n}}(\zeta) \leq -c , \quad\quad \forall \zeta \in B(\gamma_{\rho}(s_{0}), \varepsilon) .
\end{equation}

Moreover, for $\eta \in ]0, \varepsilon[$ small enough, using again the continuity of the Melrose-Sj\"ostrand flow, we may find $V_{\rho}$, a neighborhood of  $\rho$ in $T^{*}\R^{n+1}$  such that for all  $\rho' \in V_{\rho}\cap T^{*}\L \cap Char(P_{A})$,  

\begin{equation}\label{hyp-8}
\gamma_{\rho'}(s_{0}) \in B(\gamma_{\rho}(s_{0}), \eta) .
\end{equation}
In this setting, two cases may occur : 

i) $\gamma_{\rho'}(s_{0})$ is a boundary point and necessarily $r(\gamma_{\rho'}(s_{0})) \geq 0$. If $r(\gamma_{\rho'}(s_{0})) > 0$ then $\gamma_{\rho'}(s_{0})$ is a hyperbolic point. Otherwise, $r(\gamma_{\rho'}(s_{0}))=0$ and then it's a glancing strictly gliding point thanks to \eqref{hyp-7}.

ii) $\gamma_{\rho'}(s_{0})$ is an interior point (see Figure \ref{sgliding} below ).
\begin{figure}[h!]\label{sgliding}
\includegraphics[scale=0.4]{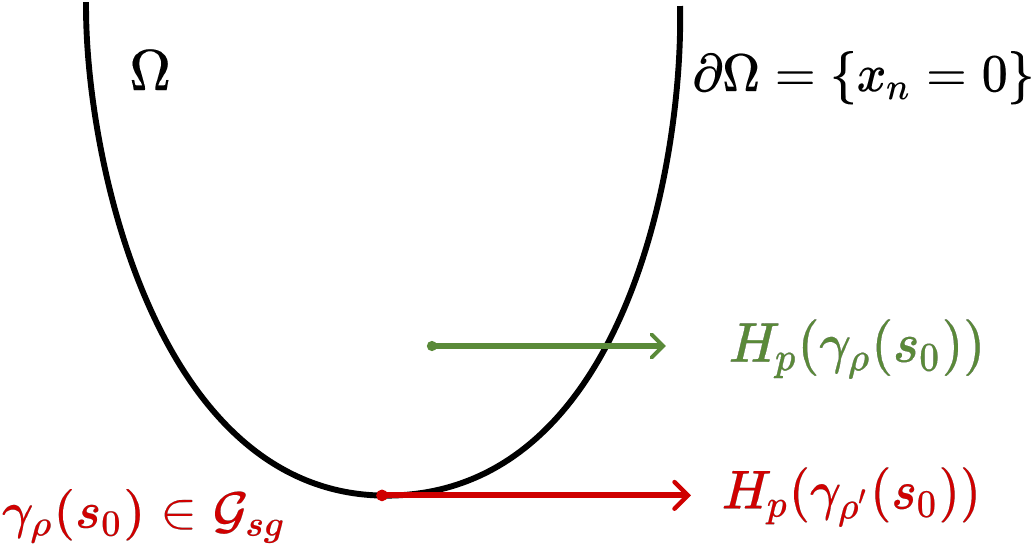}
  \caption {Strictly gliding points}
\end{figure}

 In this case,  using the Hamiltonian field $H_{p_{A}}$ , we get :
\begin{equation}\label{hyp-9}
\frac{dx_{n}}{ds}(\gamma_{\rho'}(s_{0})) = -2 \xi_{n}(\gamma_{\rho'}(s_{0})) \leq 2 \eta .
\end{equation}
Thus, if we denote in short $x_{n}(s)= x_{n}(\gamma_{\rho'}(s))$, we can perform a Taylor expansion  and get in vue of \eqref{hyp-7} :
\begin{equation}\label{hyp-10}
\left\{ 
\begin{array}{c}
x_{n}(s) = x_{n}(s_{0})+  \frac{dx_{n}}{ds}(s_{0})(s-s_{0}) + \frac{1}{2} \frac{d^{2}x_{n}}{ds^{2}}(s_{0})(s-s_{0})^{2} + o(s-s_{0})^{2} 

\\
\\
\leq \eta + 2 \eta (s-s_{0}) - c (s-s_{0})^{2} + o(s-s_{0})^{2} .
\end{array}%
\right. 
\end{equation}
Similarly, we obtain for the $\xi_{n}$ - component of $\gamma_{\rho'}(s)$ :
\begin{equation}\label{hyp-11}
\left\{ 
\begin{array}{c}
\xi_{n}(s) = \xi_{n}(s_{0})+  \frac{d\xi_{n}}{ds}(s_{0})(s-s_{0}) + o(s-s_{0})

\\
\\
\geq - \eta + c (s-s_{0})  + o(s-s_{0})
\end{array}%
\right. 
\end{equation}
From \eqref{hyp-10} we deduce that $\gamma_{\rho'}(s)$ intersects the boundary before the time $s_{1}$ such that $ s_{1}-s_{0} \approx \frac{1}{\sqrt{c}}\eta^{1/2}$. Furthermore, we conclude from \eqref{hyp-11} that $\xi_{n}(s)\geq \frac{\sqrt{c}}{2} \eta^{1/2}$ for $s$ close to $s_{1}$, which means that $\gamma_{\rho'}(s_{1})$ is a hyperbolic point of the boundary $\Gamma'$. Finally, we finish the argument by taking $\eta >0$ such that $\frac{1}{\sqrt{c}}\eta^{1/2} < T-T_{0}$.
\end{itemize}

The proof of Lemma \ref{bound-nbhood} is now complete.
\end{proof}

\subsection{First computations}\label{computation}

We consider a family of  pseudo-differential symbols in the class $\mathcal{A}^{0}$ introduced in section \ref{pdo} above, tangential and classical. Since the result we seek is of local nature, we work in a system of geodesic coordinates near  the boundary $\d\L$ and choose these symbols in the form  $q=q(x_{n},x',t,\xi',\tau)$,  and of class $C^{\infty}$ with respect to $x_{n}$, real valued, compactly supported in $(t,x',x_{n})$, and independent of $x_{n}$ in a strip $ \{\vert x_{n}\vert < \beta\}$, $\beta >0$ small enough. For instance, one may take $q$ in the form $q(x_{n},x',t,\xi',\tau)=\varphi(x_{n})\tilde{q}(x',t,\xi',\tau)$, with $\varphi \in \mathcal{C}_{0}^{\infty}(\R)$,  equal to $1$ near $x_{n}=0$. We shall denote by  $Q=Q(x_{n},x',t,D_{x',t})$ the corresponding  tangential pseudo-differential operators .

In the proofs of theorem \ref{theo}, we will make successive choices of symbols  $q$.

We recall that in the system of local geodesic coordinates, the wave equation takes the form 
\begin{equation}
\partial_{n}^{2}u+R(x_{n},x',D_{x',t})u+M_{0}(x)\partial_{n}u+M_{1}(x,\partial_{x'})u  =0 .
\end{equation}
We multiply the equation  by $Q^{2}\partial_{n}\overline{u}$ and we integrate over $\L$.
\begin{equation}\label{C1}
\left\{ 
\begin{array}{c}
I_{1}=\int_{\L}\partial_{n}^{2}u \, Q^{2}\partial_{n}\overline{u} = -\int_{\d\L} \partial_{n}u \, Q^{2}\partial_{n}\overline{u} \,d\sigma - \int_{\L}\partial_{n}u \partial_{n} \,Q^{2}\partial_{n}\overline{u} 
\\ 
\\ 
=  -\int_{\d\L} \partial_{n}u \, Q^{2}\partial_{n}\overline{u} \,d\sigma - \int_{\L}\partial_{n}u [\partial_{n} ,\,Q^{2}]\partial_{n}\overline{u} -  \int_{\L}\partial_{n}u Q^{2}\partial_{n}^{2}\overline{u}
\\
\\
 =-\int_{\d\L} \partial_{n}u \, Q^{2}\partial_{n}\overline{u} \,d\sigma - \int_{\L}\partial_{n}u [\partial_{n} ,\,Q^{2}]\partial_{n}\overline{u} -  \int_{\L}Q^{2}\partial_{n}u \partial_{n}^{2}\overline{u} + \int_{\L}(Q^{2}-Q^{*2})\partial_{n}u \partial_{n}^{2}\overline{u}
 \\
 \\
 = -\int_{\d\L} \partial_{n}u \, Q^{2}\partial_{n}\overline{u} \,d\sigma - \int_{\L}\partial_{n}u [\partial_{n} ,\,Q^{2}]\partial_{n}\overline{u} -  \int_{\L}Q^{2}\partial_{n}u \partial_{n}^{2}\overline{u}
 \\
 \\
- \int_{\L} (Q^{2}-Q^{*2})\partial_{n}uR\overline{u} - \int_{\L}M_{0}(Q^{2}-Q^{*2})\partial_{n}u \partial_{n}\overline{u} - \int_{\L}(Q^{2}-Q^{*2})\partial_{n}u M_{1}\overline{u}
 
\end{array}%
\right. 
\end{equation}

\begin{align}\label{C2}
\begin{array}{c}
I_{2}=\int_{\L} Ru \, Q^{2}\partial_{n}\overline{u} = \int_{\L} Ru \, [Q^{2},\partial_{n}]\overline{u} + \int_{\L} Ru \, \partial_{n}Q^{2}\overline{u}
\\
\\
= -\int_{\d\L} Ru \, Q^{2}\overline{u}d\sigma -\int_{\L} (\partial_{n}R)u \, Q^{2}\overline{u}  -\int_{\L} \partial_{n}u \, R^{*}Q^{2}\overline{u} - \int_{\L} Ru \, [\partial_{n} ,\,Q^{2}]\overline{u}
\\
\\
= -\int_{\d\L} Ru \, Q^{2}\overline{u}d\sigma-\int_{\L} (\partial_{n}R)u \, Q^{2}\overline{u}  -\int_{\L} \partial_{n}u \, [R^{*},Q^{2}]\overline{u} -\int_{\L} \partial_{n}u \, Q^{2}R^{*}\overline{u} - \int_{\L} Ru \, [\partial_{n} ,\,Q^{2}]\overline{u}
\\
\\
=-\int_{\d\L} Ru \, Q^{2}\overline{u}d\sigma-\int_{\L} (\partial_{n}R)u \, Q^{2}\overline{u}  -\int_{\L} \partial_{n}u \, [R^{*},Q^{2}]\overline{u} -\int_{\L} Q^{2}\partial_{n}u \, R\overline{u}
\\
\\
 -\int_{\L} (Q^{*2}-Q^{2})\partial_{n}u \, R \overline{u}-\int_{\L} \partial_{n}u \, Q^{2}(R^{*}-R)\overline{u} - \int_{\L} Ru \, [\partial_{n} ,\,Q^{2}]\overline{u}.r
\end{array}%
\end{align}
Setting $f = M_{0}(x)\partial_{n}u+M_{1}(x,\partial_{x'})u$ and summarizing all the computations above, we obtain 
\begin{equation}\label{IPP}
\int_{\d\L}\partial_{n}u Q^{2}\partial_{n} \overline{u} \,d\sigma + \int_{\d\L}Ru \,Q^{2}\,\overline{u}  \,d\sigma +\int_{\L}(\partial_{n}R)u\,Q^{2}\, \,\overline{u} =  2\Re\int_{\L}fQ^{2}\partial_{n} \overline{u} - \sum_{j=1}^{8}A_{j}.
\end{equation}
We have $\int_{\L}fQ^{2}\partial_{n} \overline{u} = \int_{\L}M_{0}\partial_{n}uQ^{2}\partial_{n} \overline{u}+\int_{\L}M_{1}uQ^{2}\partial_{n} \overline{u}$. The first term of the sum reads
\begin{align}\label{C4}
\begin{array}{c}
\int_{\L}M_{0}\partial_{n}uQ^{2}\partial_{n}\overline{u} 
= -\int_{\d\L}M_{0}u Q^{2}\partial_{n}\overline{u}\, \d\sigma -\int_{\L}(\partial_{n}M_{0})uQ^{2}\partial_{n}\overline{u}   
-\int_{\L}M_{0}u [\d_{n},Q^{2}]\partial_{n}\overline{u} -\int_{\L}M_{0}u Q^{2}\partial_{n}^{2}\overline{u}
\\
\\
= -\int_{\d\L}M_{0}u Q^{2}\partial_{n}\overline{u}\, d\sigma-\int_{\L}(\partial_{n}M_{0})uQ^{2}\partial_{n}\overline{u}   -\int_{\L}M_{0}u [\d_{n},Q^{2}]\partial_{n}\overline{u}+ \int_{\L}M_{0}u Q^{2}R\overline{u}+ \int_{\L}M_{0}u Q^{2}\overline{f}.
\end{array}
\end{align}

Finally we obtain
\begin{equation}\label{IPP-2}
\int_{\d\L}\partial_{n}u Q^{2}\partial_{n} \overline{u} \,d\sigma + \int_{\d\L}Ru \,Q^{2}\,\overline{u}  \,d\sigma +\int_{\L}u\,Q^{2}\,(\partial_{n}R) \,\overline{u} =   \sum_{j=1}^{14}A_{j}
\end{equation}

\begin{rema}  In fact, we will see later that the remaining terms $A_{j}$ for $j=1,...,14$, as described below,  do not play a role in our arguments, see Corollary \ref{restes-measures} and Lemma \ref{l.o.t}.
\end{rema} 
\begin{equation}\label{restes}
\left\{ 
\begin{array}{c}
A_{1} = \int_{\L}\partial_{n}u [\partial_{n} ,\,Q^{2}]\partial_{n}\overline{u} , \quad A_{2}= -\int_{\L}\partial_{n}u (Q^{*2}-Q^{2})R\overline{u} , \quad  A_{3} =  \int_{\L}(Q^{2}-Q^{*2})\partial_{n}u M_{0}\partial_{n}u, 
\\
 \\
 A_{4} =  \int_{\L}(Q^{2}-Q^{*2})\partial_{n}u M_{1}u , \quad A_{5}=  \int_{\L} \partial_{n}u \, [R^{*},Q^{2}]\overline{u}  , \quad  A_{6} = \int_{\L} (Q^{*2}-Q^{2})\partial_{n}u \, R \overline{u}
 \\
  \\
 A_{7} = \int_{\L} \partial_{n}u \, Q^{2}(R^{*}-R)\overline{u} , \quad A_{8}= 2\Re\int_{\L}(\partial_{n}M_{0})uQ^{2}\partial_{n}\overline{u}, \quad A_{9}= 2\Re\int_{\d\L}M_{0}u Q^{2}\partial_{n}\overline{u}\, d\sigma, 
 \\
  \\
 A_{10}= 2\Re\int_{\L}M_{0}u [\d_{n},Q^{2}]\partial_{n}\overline{u} , \quad A_{11} = - 2\Re\int_{\L}M_{0}u Q^{2}R\overline{u}, \quad A_{12} = - 2\Re\int_{\L}M_{0}u Q^{2}\overline{f}, 
 \\
  \\
 A_{13}= -2\Re\int_{\L}M_{1}uQ^{2}\partial_{n} \overline{u} , \quad A_{14}= \int_{\L} Ru \, [\partial_{n} ,\,Q^{2}]\overline{u}
\end{array}%
\right. 
\end{equation}

\section{Proof of Theorem \ref{theo}}\label{sec 5}
The proof relies on a classical strategy. We first establish a relaxed observability estimate, then we drop the compact term with the help of a unique continuation argument. 

\subsection{Relaxed observation and unique continuation}\label{sec 5.1}
\begin{proposition}\label{R-Obs}
Under assumptions A1, A2 and A3,  for every $T>T_{0}$, there exists  $c > 0$ such that for every $g \in H^{1}(\d\L)$, $supp(g) \subset \overline{\Gamma}_{M}$, the solution $u$ of  \eqref{waveequation}, satisfies the observability estimate  
\begin{equation}\label{R-estimate}
\Vert g\Vert_{H^{1}(\Gamma_{M})} \leq c \Vert \d_{n}u_{\vert\d \Omega}\Vert_{L^{2}(\Gamma'_{M+T} )} + c\Vert g\Vert_{L^{2}(\Gamma_{M})}.
\end{equation}
\end{proposition}

Also, we will need the following uniqueness result.

\begin{lemma}\label{UC}
Assume that estimate \eqref{R-estimate} holds true for all $T>T_{0}$. Then  for $g \in H^{1}(\d\L)$ with $supp(g) \subset \overline{\Gamma}_{M}$, if the solution $u$ to system \eqref{waveequation} satisfies $\d_{n}u_{\vert \d\Omega} \equiv 0 $ on $\Gamma'_{M+T}$, then $u$ vanishes identically.
In particular,  $ g \equiv 0$.
\end{lemma}

The proof of Lemma \ref{UC} is given at the end of this section and the proof of Proposition \ref{R-Obs} will be the purpose of Section \ref{sec.5.2}. Here,  we first show how we can conclude the proof of Theorem \ref{theo} using these results.

For this , we use a contradiction argument. Assume that estimate \eqref{obs-theo} is false and consider  a sequence of boundary data $(g_{k}) \in H^{1}(\d\L)$, $supp(g_{k})\subset \overline{\Gamma}_{M}$,   and $(u_{k})$ the sequence of associated solutions, with
\begin{equation}\label{contradiction1}
 \Vert \d_{n}u_{k\vert \d\Omega}\Vert_{L^{2}(\Gamma'_{M+T})}  < \frac{1}{k}\Vert g_{k}\Vert_{H^{1}(\Gamma )}.
\end{equation}
The sequence $v_{k}=\Vert g_{k}\Vert_{H^{1}(\Gamma )}^{-1}u_{k}$ then satisfies

\begin{equation}\label{contradiction2}
\left\{ 
\begin{array}{c}
P_{A} v_{k}=0, \quad v_{k\vert \Gamma'}=0, \quad \Vert v_{k\vert\d\Omega}\Vert_{H^{1}(\Gamma )} = 1,
  \text{and} \quad \Vert \d_{n}v_{k\vert \d\Omega}\Vert_{L^{2}(\Gamma'_{M+T})} < \frac{1}{k} \, .
\end{array}%
\right. 
\end{equation}

The sequence $(v_{k})$ is bounded in the energy space $C^{0}((0,M+T), H^{1}(\Omega)) \cap C^{1}((0,M+T), L^{2}(\Omega))$ accordingly to \eqref{energy}, thus we may assume that it converges weakly in the cylinder $\L_{M+T}$ to some function $v \in H^{1}(\L_{M+T})$.

In the same way, we assume that  the sequence $\tilde{g}_{k}=v_{k\vert\d\Omega}$  weakly converges to some $\tilde{g}$ in $H^{1}(\Gamma)$, with $supp(\tilde{g}) \subset \overline{\Gamma}_{M}$. Passing then to the limit $k\rightarrow \infty$ in \eqref{contradiction2}, we obtain
\begin{equation}\label{weaklimit}
P_{A} v=0, \quad v_{\vert \d\Omega}= \tilde{g},   \quad \text{and}\quad  \d_{n}v_{\vert \d\Omega} = 0 \quad \text{on} \,\, \Gamma'_{M+T}.
\end{equation}
The unique continuation result of lemma \ref{UC} then gives that the weak limits $v$ and $\tilde{g}$ vanish identically. Coming back then to Proposition \ref{R-Obs} and plugging $v_{k}$ and $\tilde{g}_{k}$ in estimate \eqref{R-estimate}, we get the contradiction 
\begin{equation*}
1 \leq  c \Vert \tilde{g}_{k}\Vert_{L^{2}(\Gamma_{M})} \longrightarrow 0 \quad \text{as} \quad k \rightarrow \infty
\end{equation*}
thanks to the compact imbedding of $H^{1}(\Gamma_{M})$ into $L^{2}(\Gamma_{M})$.

\begin{proof}[\textbf{Proof of the unique continuation}]

The proof is based on a classical argument of functional analysis. For $a\geq 0$ and $g \in H^{1}(\d\L)$ with $supp(g) \subset \overline{\Gamma}_{M}^{a} =: [-a,M]\times \overline{O}$, consider the system  
\begin{equation}
\left\{ 
\begin{array}{c}
P_{A}u=\partial _{t}^{2}u-\sum_{i,j=1}^{n}\partial_{x_{j}}(a _{ij}(x)\partial_{x_{i}}u)=0\quad \text{in } \L 
\\ 
u(t,.)=g(t,.)\quad \text{on } \d\L 
\\
u(-a,.)=\partial _{t}u(-a,.)=0 \quad \text{in } \Omega.
\end{array}%
\right.  \label{waveequation-a}
\end{equation}
Clearly, the solutions of \eqref{waveequation-a} satisfy a relaxed observability estimate similar to \eqref{R-Obs}, namely
\begin{equation}\label{R-estimate-a}
\Vert g\Vert_{H^{1}(\Gamma_{M}^{a})} \leq c \Vert \d_{n}u_{\vert \d\Omega}\Vert_{L^{2}(\Gamma'^{a}_{M+T})} + c\Vert g\Vert_{L^{2}(\Gamma_{M}^{a})}.
\end{equation}
for any $T>T_{0}$ and some $c>0$. Here we have denoted $\Gamma^{a}_{M} = (-a,M)\times O$ and $\Gamma'^{a}_{M+T} = (-a,M+T)\times O'$.

Let us introduce  the set
\begin{equation}
\mathcal{N}_{a}(T)= \Big\{ 
g \in H^{1}(\d\L), \,\,supp(g) \subset \overline{\Gamma}_{M}^{a}, \,\,  u=u(g) \,\, \text{solves} \,\, \eqref{waveequation-a} \,\,\, \text{and} \,\,  \d_{n}u_{\vert \Gamma'^{a}_{M+T}} \equiv 0 \Big\}
\end{equation}
First we notice that thanks to \eqref{b-energy}, $\mathcal{N}_{a}(T)$ is a closed subset of $H^{1}(\Gamma_{M}^{a})$. In addition, applying the relaxed observability \eqref{R-estimate-a} to  an element  of $\mathcal{N}_{a}(T)$ gives 
\begin{equation*}
\Vert g\Vert_{H^{1}(\Gamma_{M}^{a})} \leq  c \Vert g\Vert_{L^{2}(\Gamma_{M}^{a})} .
\end{equation*}
Using the compact imbedding $H^{1}(\Gamma_{M}^{a}) \hookrightarrow L^{2}(\Gamma_{M}^{a})$, this implies that $\mathcal{N}_{a}(T)$ has a finite dimension,  and thus is complete for any norm.

Now we come back to the initial problem. We pick $g \in \mathcal{N}_{0}(T)$, i.e $g \in H^{1}(\d\L)$  with support in $\overline{\Gamma}_{M}$,  and we consider $u$,  the associated solution of \eqref{waveequation}. Notice first that $g \in \mathcal{N}_{a}(T)$ for all $a>0$. In what follows, we fix $a>0$. In addition, for  $\delta = \frac{1}{2}(T-T_{0})$, we remark  that estimate \eqref{R-estimate-a} is also satisfied by all functions $h \in \mathcal{N}_{a}(T-\delta)$.    Moreover, for all $\varepsilon <\min(\delta, a)$, the function $g(t+\varepsilon, .)$ lies in $ \mathcal{N}_{a}(T-\delta)$. We also have 
\begin{equation*}
h_{\varepsilon}=\frac{1}{\varepsilon}(g(t+\varepsilon, .) -g(t, .)) \mathop{\to}_{\eps \to 0^+} \frac{\d g}{\d t} \quad \text{in} \quad L^{2}(\Gamma_{M}^{a}) .
\end{equation*}
As a consequence, the sequence $(h_{\varepsilon})_{\varepsilon>0}$ is a Cauchy sequence in $\mathcal{N}_{a}(T-\delta)$ endowed with the norm $\Vert . \Vert_{L^{2}(\Gamma_{M}^{a})}$. As all norms are equivalent , the sequence $(h_{\varepsilon})_{\varepsilon>0}$ is thus also a Cauchy sequence in $\mathcal{N}_{a}(T-\delta)$ endowed with the norm $\Vert . \Vert_{H^{1}(\Gamma_{M}^{a})}$, which yields $\frac{\d g}{\d t} \in \mathcal{N}_{a}(T-\delta)$.  In particular, $\frac{\d g}{\d t} \in H^{1}(\Gamma_{M}^{a})$. This distribution is supported in $\overline{\Gamma}_{M}$, we get therefore $\frac{\d g}{\d t} \in \mathcal{N}_{0}(T-\delta)$. Finally if $u(\frac{\d g}{\d t})$ denotes the solution of system \eqref{waveequation-a} with boundary data $\frac{\d g}{\d t}$, we write 
\begin{equation*}
\d_{n}\big(u(\frac{\d g}{\d t})\big) = \d_{n}\big(\frac{\d u(g)}{\d t}\big) = \d_{t}\big(\frac{\d u(g)}{\d n}\big) = 0 \quad \text{on} \quad (0,M+T)\times O' .
\end{equation*}

Therefore we obtain that $\frac{\d g}{\d t} \in \mathcal{N}_{0}(T)$.

To summarize, we have proved that the time derivative $\frac{\d}{\d t}$ defines a linear operator on the finite dimensional space $\mathcal{N}_{0}(T)$. But we notice that this operator has  no eigenvalue. Indeed,  for $g\in \mathcal{N}_{0}(T)$, we have $supp(g) \subset \overline{\Gamma}_{M}$; therefore for all $\lambda \in \C $, the only solution of system 
\begin{equation*}
\frac{\d g}{\d t} = \lambda g, \quad\quad g(0,.) =0
\end{equation*}
is the trivial one $g \equiv 0$. This concludes the proof of Lemma \ref{UC}. 

This also concludes the proof of Theorem \ref{theo} assuming the relaxed observation estimate \eqref{R-estimate}. Accordingly, the next section is dedicated to the proof of Proposition \ref{R-Obs}.
\end{proof}

\subsection{Proof of the relaxed observation}\label{sec.5.2}

In order to establish estimate \ref{R-Obs}, we use a contradiction argument. Assume that inequality \eqref{R-estimate} is false and consider  a sequence of boundary data $(g_{k}) \in H^{1}(\d\L)$, $supp(g_{k})\subset \overline{\Gamma}_{M}$,   and $(u_{k})$ the sequence of associated solutions, with
\begin{equation}\label{contradiction1}
 \Vert \d_{n}u_{k\vert \d\Omega}\Vert_{L^{2}(\Gamma'_{M+T})} + \Vert g_{k}\Vert_{L^{2}(\Gamma_{M} )}< \frac{1}{k}\Vert g_{k}\Vert_{H^{1}(\Gamma_{M} )}.
\end{equation}

The sequence $v_{k}=\Vert g_{k}\Vert_{H^{1}(\Gamma )}^{-1}u_{k}$ then satisfies
\begin{equation}\label{weaklimit}
P_{A} v_{k}=0, \quad v_{k\vert \d\Omega}= \Vert g_{k}\Vert_{H^{1}(\Gamma_{M} )}^{-1}g_{k},  \quad \text{and} \quad   \Vert \d_{n}v_{k\vert \d\Omega}\Vert_{L^{2}(\Gamma'_{M+T})}\rightarrow 0.
\end{equation}
$(v_{k})$ is bounded in $H^{1}(\L_{T})$ and $(v_{k\vert \d\Omega})$ is bounded in $H^{1}(\Gamma_{M})$. Therefore we may assume that $(v_{k})$ weakly converges to some $v$ in $H^{1}(\L_{T})$ and $(v_{k\vert\d\Omega})$ weakly converges to some $\tilde{g}$ in $H^{1}(\Gamma_{M})$. Equations \eqref{weaklimit} then provides
\begin{equation}\label{weaklimit-1}
P_{A} v=0, \quad v_{\vert \d\Omega}= \tilde{g},  \quad \text{and} \quad    \d_{n}v_{\vert\d\Omega} = 0,
\end{equation}
and Lemma \ref{UC}  implies that $v$ and $v_{\vert \d\Omega}= \tilde{g}$ vanish identically. Thus, the weak limits are both equal to $0$.

Our goal, will be to prove that in the contradiction setting assumed above, the sequence $(v_{k\vert \d\Omega})$  strongly converges to $0$ in $H^{1}(\Gamma)$, which is a impossible  since  $\Vert v_{k\vert \d\Omega}\Vert_{H^{1}(\Gamma_{M})} = 1$ accordingly to \eqref{weaklimit}.

For this purpose, we make use of a classical strategy.  Following Burq-Lebeau \cite{BurqLebeau}, and coming back to the notation $u_{k}$ instead of $v_{k}$, we  attach to $(u_{k})$ a microlocal defect measure  in $H^{1}(\L_{M+T})$ denoted by $\mu$. 

 Also,  we  attach to $(g_{k})$ a microlocal defect measure  on the boundary, in  $H^{1}(\d\L)$, denoted by $\tilde{\mu}$.   Finally, the sequence $\d_{n}u_{k\vert \d\Omega}$  weakly converges to $0$ in $L^{2}_{loc}(\d\L)$. So we attach to it a microlocal defect measure  in $L^{2}_{loc}(\d\L)$ denoted by $\nu$. 

Notice, that in the contradiction setting of \eqref{weaklimit}, the measure $\nu$ vanishes identically over $\Gamma'_{M+T}$.

Finally, we will prove in several steps,  that in the contradiction setting assumed above, the measure $\tilde{\mu}$ vanishes identically  on $\Gamma_{M}$. Notice that in the different intermediate results we will prove below, we    use this contradiction setting, without explicitly  referring to it.

\subsection{Properties of the measures}\label{measures-prop}

In the sequel we consider  $W$  an interior neighborhood of the boundary $\overline{\Gamma}$  as introduced in Section \ref{subsec.4}. We recall that $W = \R\times (V\cap \Omega)=(\R\times V)\cap \L$ where $V$ is an open subset of $\R^{n}$, neighborhood of the spatial boundary $O\subset \d\Omega$.  We set  
\begin{equation}\label{nbhood-0}
W^{\d}= (\R\times V)\cap \d\L .
\end{equation}
In addition, for $J$  an open interval of $\R$ such that $[0,M] \subset J$, we  denote
\begin{equation}\label{nbhood}
W_{J}=\{(t,x) \in W, \,\, t \in J\} \quad \text{and} \quad W^{\d}_{J}=\{(t,x) \in W^{\d}, \,\, t \in J\}.
\end{equation}
The neighborhood $W$ and the interval $J$ will be fixed in the next Proposition.

\begin{proposition}\label{interior}
 Under assumptions A1 and  A2, for every $T>T_{0}$, there exist $W$ and $J$ as above such that  the measure $\mu$ vanishes identically near any interior point of $W_{J}$.
\end{proposition}
\begin{proof}
Consider $T>T_{0}$. We take the  interior neighborhood  $W$ of $\Gamma$ satisfying the conclusion of  Lemma \ref{bound-nbhood} with $\frac{T+T_{0}}{2}$. In addition,  we chose $J=]-\alpha, M+\alpha[$, where $0<\alpha<\frac{T-T_{0}}{2}$. And we prove that $\rho \notin$ supp$(\mu)$ for all  $\rho \in T^{*}W_{J}$. This fact is obvious if $\rho$ is an elliptic point, thanks to the classical property of microlocal elliptic regularity. If $\rho \in Char(P_{A})$,   let $\gamma = \gamma(s)$ be the generalized half bicharacteristic starting at $\rho$ and satisfying (SGCC). We know that for some $s_{0}$ ( say  $0< s_{0}<\frac{T+T_{0}}{2}$ ), $\gamma(s_{0})=(t_{0}, x_{0}, \tau_{0},\xi_{0})$ is a strictly gliding point of the boundary $\Gamma_{M+T}'$. Consider $U_{0}$  a small neighborhood of $(t_{0} ,x_{0})$ in $\R^{n+1}$  and denote by  $\underline{u}_{k}$ the canonical extension of ${u}_{k}$ to $\R^{n+1}$, i.e  $\underline{u}_{k}=u_{k}$ in $\L$ and  $\underline{u}_{k}=0$ elsewhere. We have
\begin{equation}
\left\{ 
\begin{array}{c}
\underline{u}_{k} \rightharpoonup 0 \quad \text{in } H^{1}(U_{0}) \quad \text{weakly }
\\
\\
u_{k \vert\d\Omega}=0 \quad \text{on } U_{0}\cap\d\L  \quad \text{and } \d_{n}u_{k \vert\d\Omega} \longrightarrow 0  \quad \text{on } U_{0}\cap\d\L \quad \text{strongly }.
\end{array}%
\right.  \label{equationextension}
\end{equation}

Accordingly to the lifting lemma of Bardos, Lebeau and Rauch \cite[Theorem 2.2]{B-L-R} or Burq \cite[Lemme 2.2]{Burq}, we know that  $\underline{u}_{k}$ strongly converges to $0$ in $H^{1}$ microlocally at  $\gamma(s_{0})$. Therefore we deduce that $\gamma(s_{0}) \notin supp(\mu)$ thanks to the work of Aloui \cite[Lemme 3.1]{Aloui}. Now, accordingly to (SGCC), for $0\leq s \leq s_{0}$, the bicharacteristic $\gamma(s)$ doesn't intersect the boundary $\Gamma$. It may only intersect $\d\L\setminus\overline{\Gamma}$, on which we have homogeneous Dirichlet condition $u_{k_\vert\partial\Omega}=0$. Consequently, the measure propagation result of Lebeau \cite{Lebeau} or Burq-Lebeau \cite{BurqLebeau}  is valid.  Starting then backward from  $\gamma(s_{0})$, and using  the propagation  of the measure $\mu$, we obtain that $\rho \notin$ supp$(\mu)$. Finally, the case $s_{0}<0 , \, 0<\vert s_{0}\vert <\frac{T+T_{0}}{2}$,  can be treated in a similar way.
\end{proof}
\begin{rema}
In the rest of the proof, the neighborhood $W$ and the interval $J$ are fixed as in the proof of Proposition \ref{interior} above.
\end{rema}

\begin{proposition}\label{hyperbolic}
Under assumptions A1 and  A2, the measures $\mu, \nu$ and $\tilde{\mu}$ vanish on the hyperbolic set of the boundary $W^{\d}_{J}$.
\end{proposition}
\begin{proof}
The fact that $\mu \bold{1}_{\H}=0$ is proved in Burq-Lebeau paper ( see \cite[Lemma 2.6]{BurqLebeau} ) and is independent of the boundary condition. It only needs the weak convergence of the sequence $(u_{k})$ to $0$ in $H^{1}_{loc}(\L)$. On the other hand, since $\mu =0$ in the interior of $W_{J}$ thanks to Proposition \ref{interior},  the two hyperbolic fibers incoming to and outcoming from any hyperbolic point $\rho_{0}$ of the boundary  $W^{\d}_{J}$ are not charged, i.e they don't intersect $supp(\mu)$. Therefore, the Taylor pseudo-differential factorization ( see for instance Burq-Lebeau \cite[ Appendix]{BurqLebeau} ), shows that microlocally near $\,\rho_{0}$,  $g_{k} = u_{k\vert\d\Omega} \rightarrow 0$ in $H^{1}$ and $\partial_{n}u_{k\vert\d\Omega} \rightarrow 0$ in $L^{2}$ strongly. So as a by-product, we get that $\rho_{0}$ is not in supp $\tilde{\mu}$ neither in supp $\nu$.
\end{proof}

At this step, we can already conclude the proof of Theorem \ref{theo} under assumption A3.a.
\begin{corollary}\label{proof A3.a}
Under assumptions A1,  A2 and A3.a, the measure $\tilde{\mu}$ identically vanishes on the boundary $W^{\d}_{J}$.
\end{corollary}
\begin{proof}
This result is a byproduct of Proposition \ref{hyperbolic} and we develop it for the convenience of the reader. First we recall a classical property of micolocal defect measures, namely the microlocal elliptic regularity. 
Let $\chi=\chi(t,x',\tau,\xi')$ and  $\psi=\psi(t,x',\tau,\xi')$ two  0-order pseudo-differential symbols supported in $T^{\ast}(\d\L)_{\vert W^{\d}_{J}}\setminus Char B_{\alpha}$, such that $\chi \equiv 1$ on $supp(\psi).$  It's classical that one can find a pseudo-differential operator $B_{-\alpha}$, of order $(-\alpha)$ on $\d\L$ such that 
\begin{equation}
B_{-\alpha} B_{\alpha}\chi(t,x',D_{t},D_{x'}) = \psi(t,x',D_{t},D_{x'}) + R_{-\infty}
\end{equation}
where $R_{-\infty}$ is infinitely smoothing. Consequently, can write the elliptic estimate
\begin{equation}\label{elliptic-0}
\Vert \psi(t,x',D_{t},D_{x'})g_{k}\Vert_{H^{1}(\d\L)}  \leq  c_{0} \Vert B_{\alpha}\chi(t,x',D_{t},D_{x'})g_{k}\Vert_{H^{1-\alpha}(\d\L)} + c_{1}\Vert g_{k}\Vert_{L^{2}(\d\L)}
\end{equation}
for some constants $c_{0}, c_{1} >0$. Therefore
\begin{equation}
\Vert \psi(t,x',D_{t},D_{x'})g_{k}\Vert_{H^{1}(\d\L)} \leq c_{0} \Vert [B_{\alpha},\chi(t,x',D_{t},D_{x'})]g_{k}\Vert_{H^{1-\alpha}(\d\L)} +  c_{1}\Vert g_{k}\Vert_{L^{2}(\d\L)} \leq c_{2}\Vert g_{k}\Vert_{L^{2}(\d\L)}
\end{equation}
for some $c_{2}>0$.
We then deduce that $\psi(t,x',D_{t},D_{x'})g_{k} \rightarrow 0$ strongly in $H^{1}(\d\L)$,  which expresses that $supp(\tilde{\mu}) \subset Char B_{\alpha}$. Now , $Char B_{\alpha} \subset \H$ thanks to  assumption A3.a, and $\tilde{\mu} \equiv 0$ on $\H$ accordingly to Proposition \ref{hyperbolic}. Therefore, $\tilde{\mu}$ vanishes identically.

\end{proof}
The proof of Theorem \ref{theo} under assumption A3.a is  complete.

Let us now continue the proof of Theorem \ref{theo} under assumption A3.b. 

Denote by $A_{j}^{k}$ the terms of \eqref{restes} where we set $u_{k}$ instead of $u$, and consider a pseudo-differential symbol $q=\sigma (Q) \in \mathcal{A}^{0}$ ( see Section \ref{pdo}),  chosen as in Section \ref{computation}.   

\begin{corollary}\label{restes-measures}
Under assumptions A1 and  A2, if  $q=\sigma (Q)$ is compactly supported in $W_{J}$,  we have 
\begin{equation}\label{restes-nuls}
\lim_{k\rightarrow \infty} A_{j}^{k} = 0, \quad\quad \forall  j \in \{1, 8, 9, 10, 12, 14\} .
\end{equation}
\end{corollary}
\begin{proof}
We recall that the symbol $q=\sigma(Q)$ is  independent of $x_{n}$ in a strip $ \{\vert x_{n}\vert < \beta\}$, $\beta >0$ small. More precisely, we  take $q$ in the form $q(x_{n},x',t,\xi',\tau)=\varphi(x_{n})\tilde{q}(x',t,\xi',\tau)$, with $\varphi \in \mathcal{C}_{0}^{\infty}(\R)$,  equal to $1$ near $x_{n}=0$. Therefore, if we choose $\beta $ small enough, and assume that $\tilde{q}$ is supported in time in the interval $J$,  the symbol of the bracket operator $[\d_n,Q^{2}]$ is of order $0$ and compactly supported in the interior of $W_{J}$. Thus, $\lim_{k\rightarrow \infty} A_{j}^{k} = 0$ for $j \in \{1, 10\}$ thanks to Proposition \ref{interior}. The terms $A_{j}^{k}, \,\, j=8,9,12$ are trivial. 
\end{proof}

\begin{rema}
In the rest of the proof, we will work henceforth, with this choice of symbol $q$, and we will choose successively, the localization of its support.
\end{rema}

Now, for the convenience of the reader, we recall  the following result due to Burq-Lebeau \cite{BurqLebeau}.

In the system of geodesic coordinates introduced above, consider  the function $\theta$ defined $\mu$-almost everywhere on $S\hat{Z}$
\begin{equation}\label{measure-calculus-def}
\theta = \frac{\xi_{n}}{\vert (\tau,\xi')\vert} \,\,\, \text{in} \quad x_{n}>0 ,  \quad\quad \quad\quad  \theta = i\frac{\sqrt{-r_{0}}}{\vert (\tau,\xi')\vert} \,\,\, \text{in} \quad \E\cup\G .
\end{equation}

\begin{lemma}\cite[Lemma 2.7]{BurqLebeau} 
Let $Q_{j}\in \mathcal{A}^{j}, \, j=1,2$ be  tangential pseudo-differential operators  with principal symbols $\sigma(Q^{j})=q_{j}$. Then we have with $\lambda^{2}=\vert (\tau,\xi')\vert^{2}(1+\vert \theta\vert^{2})$ 
\begin{equation}\label{measure-calculus}
lim_{k\rightarrow \infty}\Big((Q_{2} -iQ_{1}\d_{n})u_{k} \, \vert \, u_{k}\Big)_{L^{2}(\L)} = \Big\langle \mu , \lambda^{-2}(q_{2}+q_{1}\theta\vert (\tau,\xi')\vert )\Big\rangle
\end{equation}
\label{calculus}
\end{lemma}

\begin{proposition}\label{elliptic}
The measure $\mu$ vanishes on the elliptic set of the boundary $W^{\d}_{J} $.
\end{proposition}

\begin{proof}
The elliptic microlocal regularity for measures or wave fronts is classical for elliptic interior points $\rho \in T^{\ast}W_{J}$ . In what concerns the elliptic set of the boundary, we will invoke a result of Burq-Lebeau (\cite[Lemma 2.6]{Lebeau} ), and we have to introduce some additional notations. 
\newline
In the framework above, they define a boundary measure $\mu_{\d}^{0}$ given by
\begin{equation}
\forall Q \in \mathcal{A}^{0}, \quad\quad lim_{k}\int_{\d\L}Qu_{k}\,\d_{n}\overline{u}_{k}d\sigma = \Big\langle \mu_{\d}^{0} , \sigma(Q)_{\vert x_{n}=0} \Big\rangle
\end{equation}
Moreover,
they provide the following link between the two measures $\mu$ and $\mu_{\d}^{0}$ :
\begin{equation}
\mu_{\d}^{0} = -2\frac{\vert\theta\vert^{2}}{1+\vert\theta\vert^{2}}\,\mu \bold{1}_{\vert x_{n=0}} .
\end{equation}
Therefore, we get 
$$
\mu_{\d}^{0}= \frac{2r_{0}(x'; \tau,\xi')}{\vert (\tau,\xi')\vert^{2}-r_{0}(x'; \tau,\xi')} \mu \, \bold{1}_{\vert x_{n=0}} \quad \text{on} \quad \E \cup \G 
$$
But, since $u_{k\vert \d\L}= g_{k}\rightarrow 0$ in $L^2_{loc}(\d\L)$ strongly and  $\d_{n}u_{k\vert \d\L}$ is bounded in $L^2_{loc}(\d\L)$, we easily get that $\mu_{\d}^{0}\equiv 0$. Consequently,  we obtain $\mu \equiv 0$ on $\E$,  since $r_{0}<0$ on this set.
\end{proof}.

\begin{rema}
\begin{enumerate}
\item Notice that for this proposition, we have used none of the assumptions $A_{j}$, $j=1,2,3$. We have only used the weak convergence  $g_{k} \rightharpoonup 0 $ in $H^{1}(\d\L)$ and subsequently $u_{k} \rightharpoonup 0 $ in $H^{1}(\L)$. 
 
\item One should be carefull that this proposition does not give any information about the behavior of the boundary data $g_{k}$ on $\E \cup \G$. In other words, we have not yet any information about $\tilde{\mu}\bold{1}_{\vert \E\cup G}$.

\item Up to now, we have proved that the measure $\mu$  vanishes in $T^{*}(W_{J})$ , i.e on interior points, and on the subset $\H \cup \E$ of  $T^{*}(W_{J}^{\d})$ . Therefore, $\mu$ is supported in the glancing set , that is  $\mu = \mu \bold{1}_{\G}$.
\end{enumerate}
\end{rema}

\begin{lemma}\label{l.o.t} 
Under assumptions A1 and  A2, and with a suitable choice of the pseudo-differential symbol $q=\sigma (Q)$, we have 
\begin{equation}\label{restes-nuls-1}
\lim_{k\rightarrow \infty} A_{j}^{k} = 0, \quad\quad \forall  j \in \{2, 3, 4, 5, 6, 7, 11, 13\} .
\end{equation}
Together with \eqref{restes-measures}, this implies that the right hand side of \eqref{IPP-2} tends to $0$ as $k\rightarrow \infty$.
\end{lemma} 

\begin{proof}
The proof essentially relies on the calculus Lemma \ref{calculus} . If we detail the limit \eqref{measure-calculus}, we can write accordingly to Propositions \ref{interior}, \ref{hyperbolic} and \ref{elliptic} 
\begin{equation}\label{measure-calculus-1}
\left\{ 
\begin{array}{c}
lim_{k\rightarrow \infty}\Big (Q_{2} u \, \vert \, u\Big)_{L^{2}(\L)} = \Big\langle \mu \bold{1}_{\G} , \lambda^{-2}q_{2}\Big\rangle
\\
\\
lim_{k\rightarrow \infty}\Big( -iQ_{1}\d_{n}u \, \vert \, u\Big)_{L^{2}(\L)} = \Big\langle \mu \bold{1}_{\G} , \lambda^{-2}q_{1}\theta\vert (\tau,\xi')\vert \Big\rangle = \Big\langle \mu \bold{1}_{\G} , i\lambda^{-2}q_{1}\sqrt{-r_{0}} \Big\rangle = 0
\end{array}%
\right. 
\end{equation}
since $r_{0} \equiv 0$ on the glancing set $\G$.

First, we take the pseudo-differential symbol $q=\sigma(Q)$ as in the proof of Corollary \ref{restes-measures}. With this choice, the terms  $A^{k}_{2}, A^{k}_{4}, A^{k}_{5}, A^{k}_{6},A^{k}_{7}, A^{k}_{11}$ and $A^{k}_{13}$ can be treated  with the second limit of \eqref{measure-calculus} since the pseudo-differential operator $(Q^{2}-Q^{*2})$, resp. $(R-R^{*})$ is   of order $\leq (-1), resp.  \, 1$.

On the other hand, the term $A^{k}_{11}$ tends to $0$ thanks to the first limit of \eqref{measure-calculus}. Finally, for the term $A^{k}_{3}$, we have just to notice that $\d_{n}u_{k}$ is bounded in $L^{2}_{x_{n}}(L^{2}_{t,x'})$ and converges weakly to $0$ in this space, and use again the fact that $(Q^{2}-Q^{*2})$ is   of order $\leq (-1)$.
\end{proof}

As a by-product, we have obtained the following lemma. We denote by $q=\sigma(Q)$ the symbol of the  pseudo-differential operator  $Q \in \mathcal{A}^{0}$.

\begin{corollary}
Under assumptions A1 and  A2, the measures $\mu, \tilde{\mu}$ and $\nu$  satisfy the following identity 
\begin{equation}\label{measures-equ-0}
 \Big\langle \nu, q^{2} \Big\rangle+\Big\langle \tilde{\mu}, \vert (\tau,\xi')\vert^{-2}q^{2}r_{0} \Big\rangle = - \Big\langle \mu \bold{1}_{\G}, \vert (\tau,\xi')\vert^{-2} \, q^{2}(\d_{n}r) \Big\rangle ,
\end{equation}
for all $0$-order symbol $q $, supported in $W_{J}$.
\end{corollary}

Now, we can conclude the study for the measure $\mu$.

\begin{proposition}\label{mu = 0}
The measure $\mu$ vanishes identically over $T^{*}(W_{J}^{\d})$.

In particular, $u_{k}\rightarrow 0$ strongly in $H^{1}(W_{J})$ up to the boundary.
\end{proposition}

\begin{proof}
The proof relies on a specific choice of the symbol $q$. First, we recall the notation 

$r_{0}(x',\tau, \xi')=\tau^{2}-\sum_{1\leq i,j \leq n-1}a_{ij}(x',0)\xi_{i}\xi_{j}$, see Section \ref{geo}. In addition, it's clear that in formula \eqref{measures-equ-0}, we adopt the notation $q=q_{\vert x_{n}=0}$. Let us then consider  $\tilde{q}_{0}\in \mathcal{C}_{0}^{\infty}(]0,T[,\R_{+})$ and a function $q_{0}\in \mathcal{C}_{0}^{\infty}(\R,\R)$, supported in $[-1,1]$, such that $q_{0}(s) = 1$ for $s \in [-1/2, 1/2]$. We set for $\varepsilon >0$
\begin{equation}
q_{\varepsilon}(t,x',\tau,\xi) = \tilde{q}_{0}(t)q_{0}\Big(\frac{r_{0}(x',\tau,\xi)}{\varepsilon \sum_{1\leq i,j \leq n-1}a_{ij}(x',0)\xi_{i}\xi_{j}}\Big)
\end{equation}
Plugging $q_{\varepsilon}$ into \eqref{measures-equ-0} and letting $\varepsilon \rightarrow 0^{+}$, we get by Lebesgue dominated convergence
\begin{equation}
 \Big\langle \nu, \bold{1}_{\G} \Big\rangle = - \Big\langle \mu \bold{1}_{\G}, \vert (\tau,\xi')\vert^{-2} \,(\d_{n}r)\Big\rangle 
\end{equation}
All  points of the glancing  set $\G=\G_{d}$ are strictly diffractive ( see \eqref{diffractive}) which gives $\d_{n}r_{\vert\G} >0$. Therefore  the two members of this identity are of opposite sign and thus  both are equal to zero. Consequently, the measure $\mu$ vanishes identically.
\end{proof}

\begin{rema}
\begin{enumerate}
\item Finally, summarizing previous results, we obtain that the measures equation \eqref{measures-equ-0} reads as follows :
\begin{equation}\label{measures-equ-1}
 \Big\langle \nu \bold{1}_{\E\cup\G}, q^{2} \Big\rangle+\Big\langle \tilde{\mu}\bold{1}_{\E\cup\G}, \vert (\tau,\xi')\vert^{-2}q^{2}r_{0} \Big\rangle = 0
\end{equation}
for all 0-order symbol $q $ , supported in $W_{J}$.

\item Roughly speaking, this formula tells us that we have two ways to prove that $\tilde{\mu} \equiv 0$. Either, we set a condition on the data $g$ itself, in other words, we make use of assumption A3.a or A3.b, or we we use a condition linking the two boundary data $\d_{n}u_{\vert\d\L}$ and $u_{\vert\d\L}=g$, which is assumption A3.c. 

\end{enumerate}

\end{rema}

\subsection{End of the proof of Theorem \ref{theo}}
Here we have reached the point where, for the first time,  we make use of assumptions  A3.b or A3.c .

\begin{proposition}
 Under assumptions A1, A2 and A3.b, the measures  $\tilde{\mu}$ and $\nu$ vanish identically  on  the set $\E\cup \G$ and hence on the boundary $\d\L$.
\end{proposition}
\begin{proof}
In the setting of assumption A3.b,  for  every $t \in J$ we can write the classical elliptic estimate
\begin{equation}\label{elliptic_{2}}
\Vert g_{k}(t,.)\Vert_{H^{1}(\d\Omega)}  \leq  c_{0} \Vert c(t,x',D_{x'})g_{k}(t,.)\Vert_{H^{1-\alpha}(\d\Omega)} + c_{1}\Vert g_{k}(t,.)\Vert_{L^{2}(\d\Omega)}=c_{1}\Vert g_{k}(t,.)\Vert_{L^{2}(\d\Omega)}
\end{equation}
for some constants $c_{0}, c_{1} >0$ independent of $t \in J$. We deduce that uniformly with respect to $t \in J$,
$$
\Vert D_{x'_{j}}g_{k}(t,.)\Vert_{L^{2}(\d\Omega)} \rightarrow 0 \quad \text{for} \quad k \rightarrow \infty
$$
Therefore, integrating on $t$ and taking the limit $k \rightarrow \infty$, we can write 
\begin{equation}\label{glancing2}
\Big\langle \tilde{\mu}, \vert(\tau,\xi')\vert^{-2}\vert\xi'\vert^{2} \Big\rangle =0
\end{equation}
and this yields 
\begin{equation}\label{glancing}
\Big\langle \tilde{\mu}, \vert(\tau,\xi')\vert^{-2}q^{2}\tau^{2} \Big\rangle = \Big\langle \tilde{\mu}\bold{1}_{\E\cup\G}, \vert(\tau,\xi')\vert^{-2}q^{2}\tau^{2} \Big\rangle = 0
\end{equation}
since $\tau^{2}\leq c\vert\xi'\vert^{2}$ in $\E\cup\G$. Together  with the result of Proposition \ref{hyperbolic}, this gives $\tilde{\mu} \equiv 0$ and $\nu \equiv 0$ accordingly to \eqref{measures-equ-1}.

This completes the proof of Theorem \ref{theo} under assumption A3.b.
\end{proof}

\begin{proposition}
 Under assumptions A1, A2 and A3.c, the measures  $\tilde{\mu}$ and $\nu$ vanish identically  on  the set $\E\cup \G$ and hence on the boundary $\d\L$.
\end{proposition}
\begin{proof}
All identities we will handle in this proof take place on the boundary $\d\L$. Therefore, we will simply write $\d_{n}u_{k}$ (resp. $u_{k}$) instead of $\d_{n}u_{k\vert \d\L}$ (resp. $u_{k\vert \d\L}$). In addition, without loss of generality, we may assume that $  \mathcal{U}_{M}\subset W_{J}$. Denote $F_{k}=\d_{n}u_{k} + \d_{t}u_{k}$. Clearly, $F_{k}\rightharpoonup 0$ weakly in $L^{2}(\d\L)$. In addition, thanks to condition A3.c, $F_{k}$ is bounded in $H^{\alpha}(\mathcal{U}_{M})$, with $\alpha >0$. Therefore we may assume that 
\begin{equation}\label{measures-equ-5}
 \d_{n}u_{k} + \d_{t}u_{k} = F_{k}\rightarrow 0  \quad \text{strongly  in} \quad L^{2}(\mathcal{U}_{M}).
\end{equation}

Consider an elliptic point $\rho_{0} \in T^{*}(\mathcal{U}_{M})$. A classical analysis at elliptic points of the boundary , see for instance \cite[Appendix]{BurqLebeau}, shows that microlocally near $\rho_{0}$, we have 
\begin{equation}\label{measures-equ-6}
\d_{n}u_{k} - Op(\sqrt{-r_{0}(x',t,\tau,\xi')})u_{k} =o(1) \quad \text{in} \quad H^{1/2} ,  \quad \text{for} \quad k \rightarrow \infty .
\end{equation}
Together with \eqref{measures-equ-5}, this yields
\begin{equation}\label{measures-equ-7}
\d_{t}u_{k} + Op(\sqrt{-r_{0}(x',t,\tau,\xi')})u_{k} =o(1) \quad \text{in} \quad L^{2} ,   \quad \text{for} \quad k \rightarrow \infty .
\end{equation}
Therefore  $u_{k\vert\d\L}=g_{k}\rightarrow 0 $ strongly  in $H^{1}$ near $\rho_{0}$ since the symbol $i\tau +\sqrt{-r_{0}(x',t,\tau,\xi')}$ is elliptic near this point. Consequently  $\rho_{0} \notin supp(\tilde{\mu})$ and using  \eqref{measures-equ-6}, $\rho_{0} \notin supp(\nu)$ . Thus we obtained
\begin{equation}\label{measures-equ-7-bis}
\nu = \nu \bold{1}_{\G} \quad \text{and} \quad \tilde{\mu}= \tilde{\mu}\bold{1}_{\G} .
\end{equation}
Coming back to \eqref{measures-equ-1} and using a test symbol $q$ elliptic near $\G$, we then get $\nu \equiv 0$ since $r_{0}=0$ on $\G$.

On the other hand, if  $Q$ is a  0-order polyhomogeneous pseudo-differential operator  on $\d\L$, with symbol $q $,  real valued and supported in $\mathcal{U}_{M}$, we have
\begin{equation}
\Big(Q{^{2}} \d_{n}u_{k} \, \vert \, \d_{n}u_{k}\Big)_{L^{2}(\mathcal{U}_{M})}=\Big(Q{^{2}} \d_{t}u_{k} \, \vert \, \d_{t}u_{k}\Big)_{L^{2}(\mathcal{U}_{M})} +\Big(Q{^{2}} F_{k} \, \vert \, F_{k}\Big)_{L^{2}(\mathcal{U}_{M})} -2Re\Big(Q{^{2}} F_{k} \, \vert \, \d_{t}u_{k}\Big)_{L^{2}(\mathcal{U}_{M})}
\end{equation}
Passing to the limit in $k$ and taking into account  \eqref{measures-equ-5} and \eqref{measures-equ-7-bis}, we obtain
\begin{equation}\label{measures-equ-7}
 \Big\langle \tilde{\mu}\bold{1}_{\G}, \vert (\tau,\xi')\vert^{-2}q^{2}\tau^{2} \Big\rangle = 0
\end{equation}
for all symbol $q $. And this gives $ \tilde{\mu} \equiv 0$ since $\tau \neq 0$ near $\G$.
%

This completes the proof of Theorem \ref{theo} under assumption A3.c.
\end{proof}
   
\section{Proof of Theorem \ref{Obs-Loss-1} and Corollary \ref{Obs-Loss-3}}\label{sec 6}

The proof is based on the wave front propagation theorem of Melrose-Sj\"ostrand , see \cite{MeSj}. We start  with a general remark about solutions of system \eqref{waveequation}. Consider $g \in H^{1}(\d\L)$, with support in $\Gamma_{M}=(0,M)\times O$ and assume  in addition that $WF(g)$, the $\mathcal{C}^{\infty}$-wave front of $g$, is contained in the elliptic set $\E$. First, we recall that the corresponding solution $u$ vanishes identically for $t\leq 0$.   Therefore  $u$ is of class $\mathcal{C}^{\infty}$ up to the boundary $\d\L$, outside $\overline{\Gamma}_{M}$. Indeed, consider $\rho \in T^{*}_{b}(\L)$, $\rho \notin T^{*}(\Gamma_{M})$, and denote $\gamma_{\rho}$ the generalized bicharacteristic curve issued from  $\rho$. Following this curve backward in time, one enters in the region $\{t<0\}$, say at some point $\gamma_{\rho}(-t_{0}),\, t_{0}>0$,  where $u$ is smooth. Accordingly to the description of a generalized bicharacteristic curve given in Section \ref{boundary-geometry}, we have for  $s_{0} \in [-t_{0},0]$ 

\begin{itemize} 
\item $\gamma_{\rho}(s_{0})$ is an interior point, i.e it lies in the characteristic set $ Char(P_{A})\cap T^{*}(\L)$ ,
\item $\gamma_{\rho}$ hits the boundary at a hyperbolic point for $s=s_{0}$ ,
\item $\gamma_{\rho}(s_{0})$ is a glancing point, i.e $\gamma_{\rho} \in \G$ .
\end{itemize}

In all cases, $\gamma_{\rho}(s)$  never intersects the closed set $WF(g) \subset \E$. Hence by regularity propagation (see \cite{MeSj}), $\rho \notin WF(u)$. Moreover, this propagation property yields that the $H^{\alpha}$ norm of $u$ is microlocally bounded near $\rho$, for every $\alpha \geq 1$. 

In the sequel we use this property to prove that estimate \eqref{obs-theo} fails in general.

Take $s<0$, $ \alpha >1$ and consider $F$  a closed conical subset of  $T^{*}(\Gamma_{M})$, $F\subset  \E$, and $V_{F}$ a  conical neighborhood of $F$ in $ T^{*}(\Gamma_{M}) \cap \E$. Finally, consider a symbol $b(t,x', \tau,\xi') \in \Psi_{phg}^{0}(\d\L)$, supported in $V_{F}$ and equal to $1$ on $F$. Denoting $B=b(t,x', D_{t},D_{x'})$ the corresponding pseudo-differential operator, it's classical that one can construct a sequence  of  smooth functions  $(f_{k}) \subset H^{1}(\d\L)$, compactly supported in $\Gamma_{M}$,  satisfying  
\begin{equation}\label{WF1}
  \Vert f_{k}\Vert_{H^{s}}=1 \quad \text{and} \quad f_{k}\rightharpoonup 0  \quad \text{weakly in} \quad H^{s}(\Gamma_{M}),
\end{equation} 
and 
\begin{equation}\label{WF2}
  \Vert Bf_{k}\Vert_{H^{s}} \rightarrow 1  \quad \text{for} \quad k \rightarrow \infty .
\end{equation} 
This simply means that the lack of compactness of $(f_{k})$  is located in $supp(b) \subset V_{F}\subset \E$.

We claim that  with a suitable choice of the symbol $b(t,x', \tau,\xi')$, the sequence $g_{k}= B f_{k}$ will be the key of our counterexample.

 First, consider a symbol $q \in \Psi_{phg}^{0}(\d\L)$ such that $supp (q) \cap V_{F} = \emptyset$.  Since  the composition  operator $Op(q)B$ is infinitely smoothing, it's classical that $\Vert Op(q)g_{k}\Vert _{H^{\alpha}}$ is uniformly bounded. More precisely, we have for some constant $c>0$
\begin{equation}\label{WF4}
\Vert Op(q)g_{k}\Vert _{H^{\alpha}}=\Vert Op(q)Bf_{k}\Vert _{H^{\alpha}} \leq c \Vert f_{k}\Vert _{H^{s}} = c .
\end{equation}
Moreover, accordingly to \eqref{WF1}, we obtain that $Op(q) g_{k} \rightarrow  0$ strongly in $H^{\alpha'}(\Gamma_{M})$ for all $\alpha'<\alpha$.

Let us now analyze the sequence  $(u_{k})$ of solutions to the wave system \eqref{waveequation}  with $(g_{k})$ as boundary data . 
\begin{equation}\label{elliptic-c-exple}
\left\{ 
\begin{array}{c}
P_{A}u_{k} = 0 \quad \text{in} \quad \L, \quad u_{k\vert\d\L} = g_{k}= B f_{k}
\\
\\
u_{k}(0)=\d_{t}u_{k}(0)= 0 . 

\end{array}%
\right. 
\end{equation}
%

We will  need the following Lemma. 
\begin{lemma}
Consider $s<0$ and for $c>0$ denote $\E_{c} =\{(\tau,\xi) \in \R^{n} , \,\, \vert\tau\vert \leq c\vert\xi\vert\}$. Then on the  space $\{ h \in H^{s}(\R^{n}), supp(\hat{h})\subset \E_{c} \}$,  $\Vert . \Vert_{L^{2}(\R;H^{s}(\R^{n-1}))}$ is a norm, equivalent to its natural norm $\Vert . \Vert_{H^{s}(\R^{n})}$. 
\end{lemma}
The proof is straightforward and left to the reader.

Now we choose the $0$-order pseudo-differential  operator $B$ introduced above in the form 
$$
B=b(t,x', D_{t},D_{x'}) = b_{1}(t,x')b_{2}(D_{t},D_{x'}) ,
$$
with $b_{1} \in \mathcal{C}^{\infty}_{0}(\Gamma_{M})$, $b\equiv 1$ on $F$, and $supp (b)\subset V_{F} \subset \E$. Necessarily, we have  for some $c>0$ ,

$$
supp (b_{2}(\tau,\xi'))\subset \{(\tau,\xi') \in \R^{n} , \,\, \vert\tau\vert \leq c\vert\xi'\vert\} .
$$

We then deduce that the sequence $(Bf_{k})$ is bounded in $L^{2}(0,M+T;H^{s}(O))$ . Therefore system \eqref{elliptic-c-exple} is well posed and $(u_{k})$ is bounded in $\mathcal{C}^{0}(0,M+T;H^{s}(\Omega))$ (see \cite[Th.2.7]{La-Lio-Triggiani}), 
 and thus in $H^{s}(\L_{M+T})$. Using the propagation argument developed in the beginning of this section, we see that $(u_{k})$ is bounded in $H^{\alpha}(\L_{M+T})$ up to the boundary, except on the subset $V_{F} \subset \E $. 
In particular, this sequence is bounded  in $H^{\alpha}(\mathcal{U})$ for any $\mathcal{U}$ interior neighborhood of the boundary observation region $\Gamma_{M+T}'=(0,M+T)\times O'$, ie : 
\begin{equation}\label{WF5}
\Vert u_{k}\Vert _{H^{\alpha}(\mathcal{U})} \leq c  \quad \text{for some} \quad  c>0. 
\end{equation}
Finally, since $u_{k} \rightharpoonup 0$ weakly in $H^{s}(\L)$ thanks to \eqref{WF1}, we obtain that $u_{k} \rightarrow  0$ strongly in $H^{\alpha'}(\mathcal{U})$ for any $\alpha'\in [1, \alpha[$, and this gives 
$$
\Vert \d_{n}u_{k\vert\d\Omega}\Vert _{L^{2}(\Gamma'_{M+T})} \rightarrow  0  .
$$
This concludes the proof of Theorem \ref{Obs-Loss-1}.

We come now to the proof of Corollary \ref{Obs-Loss-3}.

To this end, consider a glancing point $\omega_{0} \in T^{*}(\Gamma_{M})\cap \G$ and a sequence of elliptic points $(\omega_{k})_{k\geq 1} \subset T^{*}(\Gamma_{M})\cap \E$,  converging to  $\omega_{0}$ for $k \rightarrow \infty$ . For every $k \geq 1$, we  apply Theorem \ref{Obs-Loss-1} above with $F=\{\omega_{k}\}$. We then have a sequence $(g_{k}^{p})_{p\geq 1} \subset H^{s}(\d\L)$ supported in $\Gamma_{M}$ , weakly converging to $0$  in $H^{s}(\Gamma_{M})$,  such that
\begin{equation}\label{WF3}
  \lim_{p\rightarrow \infty}\Vert g_{k}^{p}\Vert_{H^{s}} = 1  \quad \text{and} \quad \lim_{p\rightarrow \infty}\Vert \d_{n}u^{p}_{k\vert\d\Omega}\Vert _{L^{2}(\Gamma'_{M+T})} =  0 ,  
\end{equation} 
where $(u^{p}_{k})_{p\geq 1}$ denotes the sequence of associated solutions .

For $k=1$, there exists $p_{1} \geq 1$  such that 

\begin{equation*}
 \vert \Vert g_{1}^{p}\Vert_{H^{s}} - 1\vert \leq 1/2   \quad \text{and} \quad \Vert \d_{n}u^{p}_{1\vert\d\Omega}\Vert _{L^{2}(\Gamma'_{M+T})} \leq 1/2 ,  \quad  \forall  p \geq p_{1} .
\end{equation*} 

Also, for $k=2$, there exists $p_{2} \geq p_{1}+1$  such that 
\begin{equation*}
 \vert \Vert g_{2}^{p}\Vert_{H^{s}} - 1\vert \leq 1/3   \quad \text{and} \quad \Vert \d_{n}u^{p}_{2\vert\d\Omega}\Vert _{L^{2}(\Gamma'_{M+T})} \leq 1/3 ,  \quad  \forall  p \geq p_{2} .
\end{equation*} 
And more generally, for every $k\geq 2$, there exists $p_{k} \geq p_{k-1}+1$  such that 
\begin{equation*}
 \vert \Vert g_{k}^{p}\Vert_{H^{s}} - 1\vert \leq \frac{1}{k+1}   \quad \text{and} \quad \Vert \d_{n}u^{p}_{k \vert\d\Omega}\Vert _{L^{2}(\Gamma'_{M+T})} \leq \frac{1}{k+1} ,  \quad  \forall  p \geq p_{k} .
\end{equation*} 

It's then easy to see that the sequence $(h_{k})_{k\geq 1} = (g_{k}^{p_{k}})_{k\geq 1}$ provides the counterexample of Corollary \ref{Obs-Loss-3}. Indeed, 
\begin{equation*}
  \Vert h_{k}\Vert_{H^{s}} \rightarrow 1  \quad \text{and} \quad \Vert \d_{n}u^{k}_{k\vert\d\Omega}\Vert _{L^{2}(\Gamma'_{M+T})} \rightarrow  0 ,  \quad \text{for} \quad k \rightarrow \infty .
\end{equation*} 
and the lack of compactness is contained in the limit set $\{\omega_{0}\}\subset \G$.

\end{document}